\documentclass[amscd,amssymb,verbatim,11pt]{amsart}
\numberwithin{equation}{section}
\usepackage{epsfig}
\usepackage{multido}
\usepackage{color}
\usepackage{pstricks}
\usepackage{pst-all}
\usepackage{pst-math}
\usepackage{pst-func}
\usepackage{pst-3dplot}
\usepackage{ulem}

\newtheorem{thm}{Theorem}[section]
\newtheorem{cor}[thm]{Corollary}
\newtheorem{lem}[thm]{Lemma}

\theoremstyle{definition}

\newtheorem{rem}[thm]{Remark}

\newif\ifShowLabels
\ShowLabelstrue
\newdimen\theight
\def\TeXref#1{%by M. Doob
     \leavevmode\vadjust{\setbox0=\hbox{{\tt
            \quad\quad  {\small  \bf #1}}}%
     \theight=\ht0
     \advance\theight  by  \dp0
     \advance\theight  by  \lineskip
     \kern -\theight \vbox  to
     \theight{\rightline{\rlap{\box0}}%
      \vss}%
      }}%

    {\begin{thm}\label{#1} \ifShowLabels \TeXref{#1} \fi}%
    {\end{thm}}

    {\begin{def}\label{#1} \ifShowLabels \TeXref{#1} \fi}%
    {\end{def}}

    {\begin{lem}\label{#1} \ifShowLabels \TeXref{#1} \fi}%
    {\end{lem}}

    {\begin{cor}\label{#1} \ifShowLabels \TeXref{#1} \fi}%
    {\end{cor}}

\newcommand{\eqRef}[1]%
     {\ifShowLabels \TeXref{#1} \fi
      \begin{equation}\label{#1} }

\ShowLabelsfalse

\newcommand{\vsp}{\vskip 1em}

\newcommand{\NI}{\noindent}
\newcommand{\bea}{\begin{eqnarray}}
\newcommand{\eea}{\end{eqnarray}}
\newcommand{\IR}{I\!\!R}
\newcommand{\bas}{\begin{align*}}
\newcommand{\eas}{\end{align*}}
\newcommand{\ba}{\begin{align}}
\newcommand{\ea}{\end{align}}

\newcommand{\be}{\begin{equation}}
\newcommand{\ee}{\end{equation}}
\newcommand{\ben}{\begin{eqnarray*}}
\newcommand{\een}{\end{eqnarray*}}
\newcommand{\lam}{\lambda}
\newcommand{\Om}{\Omega}

\newcommand{\tht}{\theta}
\newcommand{\p}{\partial}
\newcommand{\al}{\alpha}
\newcommand{\bt}{\beta}

\newcommand{\vep}{\varepsilon}
\newcommand{\dl}{\delta}
\newcommand{\D}{\Delta}

\newcommand{\s}{\sigma}
\newcommand{\su}{\subseteq}

\newcommand{\lamo}{\lambda_{\Om}}
\newcommand{\rh}{\rho}

\newcommand{\R}{\mathbb{R}}

\newcommand{\Df}{\D_\infty}

\title[{\large $\infty$}-Laplacian]{
An Eigenvalue problem for the Infinity-Laplacian}

\author[T. Bhattacharya and L. Marazzi]{Tilak Bhattacharya and Leonardo Marazzi}

\begin{document}

\maketitle

\begin{abstract}  In this work, we study an eigenvalue problem for the infinity-Laplacian on bounded domains. 
We prove the existence of the principal eigenvalue and a corresponding positive eigenfunction. The work also contains existence results when the parameter, in the equation, is less than the first eigenvalue. A comparison principle applicable to these problems is also proven. Some additional results are shown, in particular, that on star-shaped domains and on 
$C^2$ domains higher eigenfunctions change sign. When the domain is a ball, we prove that the first eigenfunction has one sign, radial principal eigenfunction exist and are unique up to scalar multiplication, and that there are infinitely many eigenvalues.
\end{abstract}

\section{\bf Introduction}
 
 In this work, we study a version of the eigenvalue problem  for the infinity-Laplacian on bounded domains. In a sense, this is a follow-up of the works in
 \cite {BMO1, BMO2} that discuss Dirichlet problems involving right hand sides that depend on the solution. 
 
In order to describe the problem better, we introduce some notations. Let $\Om\subset \IR^n,\;n\ge 2,$ be a bounded domain, $\overline{\Om}$ its closure and $\p\Om$ its boundary. We take $a\in C(\Om)\cap L^{\infty}(\Om),\;a>0$. We seek a pair $(\lam, u)$, $\lam$ real, and $u\in C(\overline{\Om})$ which solves
\eqRef{egpb}
\Df u+\lam a(x) u^3=0,\;\;\mbox{in $\Om$ and $u=0$ on $\p\Om$.}
\ee
We refer to $\lam$ as an eigenvalue of (\ref{egpb}) and to $u$ as an eigenfunction corresponding to $\lam$. The operator $\Df$ is the infinity-Laplacian and it is defined as 
$$\Df u=\sum_{i,j=1}^n \frac{\partial u}{\partial x_i}\frac{\partial u}{\partial x_j} \frac{\partial^2 u}{\partial x_i\partial x_j}.$$ 
Since $u$ is only continuous in $\overline{\Om}$ and the infinity-Laplacian is a nonlinear-degenerate elliptic operator,  solutions are to be understood in the viscosity sense. Questions involving the infinity-Laplacian have been attracting considerable attention recently. In particular, existence, uniqueness and local regularity have become topics of great interest. For greater motivation and context, we direct the reader to the works \cite{AMJ,BDM,C,CEG, PSSW}. Our current work is more along the lines of \cite{BMO1,BMO2, LW1, LW2}. From hereon, we will often refer to (\ref{egpb}) as the eigenvalue problem.

One of the main tasks is to be able to characterize the principal or the first eigenvalue of (\ref{egpb}). The seminal work \cite{BNV} provides us with an approach to achieving this goal.
While \cite{BNV} treats the case of the Laplacian, the ideas employed in it are general enough to be applicable to nonlinear operators, as shown in \cite{BD}. The work that comes closest to ours is in \cite{J}, which treats the case of the one-homogeneous infinity-Laplacian. One of the major discussion in \cite{BNV, BD, J} is the maximum principle when
the parameter $\lam$ is less than the first eigenvalue. Our work also addresses this issue in the context of (\ref{egpb}) and we prove analogues of some of the results
known for elliptic operators.

We also mention that there is great interest in studying the equation that arises when one takes the limit, as $p$ tends to infinity, of the first eigenvalue problem for the $p$-Laplacian. The resulting problem is often referred to as the infinity-eigenvalue problem, see for instance \cite{BK, JLM, JL}. The results in this current work, however, bear no relation to the
questions that arise from this problem. 

We have divided our work as follows. Section 2 contains preliminary results and estimates that will be needed for the existence of the first eigenvalue and a positive eigenfunction. We also prove a comparison principle that will be used quite frequently in this work. We also show that if $\lam$ is large enough then solutions to the problem in (\ref{egpb}) change sign. A related result appears in Section 5. Section 3 contains results for the case $\lam<\lamo$, where $\lamo$ stands for the first eigenvalue in (\ref{egpb}). We  prove a version of the maximum principle and show the existence of solutions to (\ref{egpb}) with non-trivial boundary data and right hand side. Section 4 contains a proof of the existence of the first eigenvalue and a corresponding positive eigenfunction. Also included here, is a result about the monotonicity of the first eigenvalues of the level sets of a positive eigenfunction on $\Om$. In Section 5, we study (\ref{egpb}) on $C^2$ domains and prove some results. This also contains a brief discussion for star-shaped domains. In particular, we show that eigenfunctions, corresponding to higher eigenvalues, change sign. It is not clear to us, at this time, if the above result holds in general domains. Also, we have been unable to decide if, in general, a first eigenfunction has one sign and if $\lamo$ is simple. A partial result appears in Section 6.
In Section 6, we take up the case of the ball and study the radial first eigenfunction when $a(x)$ is radial. Next, we discuss the radial version of the eigenvalue problem when $a(x)$ is a constant function. In particular, we prove that there are infinitely many eigenvalues that support radial eigenfunctions. In addition, we present a proof that the first eigenfunction, on the ball, has one sign and the radial first eigenfunctions are unique up to scalar multiplication.  

We thank the anonymous referee for reading the paper carefully and for his/her comments. We also thank Ahmed Mohammed for some discussions at the initial stages of this work. 
\vsp
\section{\bf Comparison principles and some preliminary estimates}

This section contains a version of a comparison principle which will be used throughout this work. We also list some estimates which will assist us in proving the existence of a first eigenvalue of (\ref{egpb}). In particular, we provide conditions under which solutions to (\ref{egpb}) with positive boundary data may have a priori bounds. As pointed out in the introduction, we also prove that solutions to (\ref{egpb}) change sign if $\lam$ is large enough.

We start with some notations. We work in $\IR^n,\;n\ge 2$, and if $x\in \IR^n$, we will sometimes write $x=(x_1,x_2,\cdots,x_n)$. By $e_1,e_2,\cdots, e_n$ we denote the unit vectors along the positive $x_1,x_2,\cdots,x_n$ axes respectively. We will use
$o$ to denote the origin. By $B_s(p),\;s>0$, we denote the ball of radius $s$ centered at $p$.  
We reserve $\lam$ to be a real number and it will represent the parameter in the differential equation in (\ref{egpb}). By $\lamo$, we will mean the first or the principal eigenvalue of the problem on the domain $\Om$. A careful definition of $\lamo$ will be provided later in Sections 3 and 4. Unless otherwise mentioned, the functions we encounter in this work will all be continuous. Also, throughout this work, all differential equations and inequalities are to be understood in the sense of viscosity, see \cite{CIL}. 

We recall that the in-ball of a domain $\Om$ is the largest ball that is contained in $\Om$, and the out-ball of $\Om$ is the smallest ball that contains $\Om$.

Let $\Om\subset \IR^n$ be a domain, $f\in C(\Om\times \IR,\IR)$ and $b\in C(\p\Om)$.
A function $u\in C(\Om)$ is said to be a viscosity sub-solution to $\Df u=f(x,u)$ or said to solve $\Df u\ge f(x,u)$, in $\Om$, if the following holds. For any $\psi\in C^2(\Om)$ such that  $u-\psi$ has a local maximum at a point $p\in \Om$, we have
$$\Df \psi(p)\ge f(p,u(p)).$$
Similarly, $u\in C(\Om)$ is said to be a viscosity super-solution to $\Df u=f(x,u)$ or said to solve $\Df u\le f(x,u)$, in $\Om$, if, for any $\psi\in C^2(\Om)$ such that $u-\psi$ has a local minimum at $q\in \Om$, we have
$$\Df \psi(q)\le f(q,u(q)).$$
A function $u\in C(\Om)$ is a viscosity solution to $\Df u=f(x,u)$, if it is both a sub-solution and a super-solution. 

We now introduce the following definitions in  relation to the problem (\ref{egpb}). We define $u\in C(\overline{\Om})$ to be a sub-solution to the problem
\eqRef{visc}
\Df u=f(x,u(x)),\;\;x\in \Om,\;\;\mbox{and $u=b$ on $\p\Om$,}
\ee
if $u$  satisfies $\Df u\ge f(x,u)$, in $\Om$, and $u\le b$ on $\p\Om$. Similarly, $u\in C(\overline{\Om})$ is a super-solution to (\ref{visc}) 
if $u$  satisfies $\Df u\le f(x,u)$, in $\Om$, and 
$u\ge b$ on $\p\Om$. We define $u\in C(\overline{\Om})$ to be a solution to (\ref{visc}), if it is both a sub-solution and a super-solution to (\ref{visc}).

Let us also note that the operator $\Df$ is reflection, rotation and translation invariant. We will also have the need to employ the radial version of $\Df u$, see Section 6. Suppose that
 for some $p\in \IR^n$ and for some $\rho>0,$ we have $u:\;B_{\rho}(p)\rightarrow \IR$. If $u(x)=u(r),$ where $r=|x-p|$, then we obtain by a differentiation that
\eqRef{radl}
\Df u =\left(\frac{du}{dr}\right)^2\frac{d^2u}{dr^2},\;\;r<\rho.
\ee
Setting $\s=3^{4/3}/4$, we also note that if  $u(x)=\s|x-p|^{4/3}$, then $\Df u=1,\;x\in B_{\rh}(p),$ in the sense of viscosity. 

We now gather various preliminary results we will need in the rest of this work. We start with a comparison principle. This is a variant of a result proven 
in \cite{BMO1}, see Lemma 4.1 therein. We provide details of the proof of this version. 

\begin{lem}\label{cmp0} Let $\Om\subset \IR^n,\;n\ge 2,$ be a bounded domain, $f:\Om\times \IR\rightarrow \IR$ and $g:\Om\times \IR\rightarrow \IR$ be continuous. 
Suppose that $u\in C(\overline{\Om})$ and $v\in C(\overline{\Om})$. \\
\NI (a) If $\sup_{\Om}(u-v)>\sup_{\p\Om}(u-v)$ and the following hold, in the sense of viscosity,
$$\Df u+f(x,u(x))\ge 0\;\;\;\mbox{and}\;\;\;\Df v+g(x,v(x))\le 0,\;\;\;\forall\;x\in \Om,$$
then there is a point $p\in \Om$ such that
$$(u-v)(p)=\sup_{\Om}(u-v)\;\;\;\mbox{and}\;\;\;\;g(p,v(p))\le f(p, u(p)).$$

\NI (b) Analogously, if $\inf_{\Om}(u-v)<\inf_{\p\Om}(u-v)$ and 
$$\Df u+f(x,u(x))\le 0\;\;\;\mbox{and}\;\;\;\Df v+g(x,v(x))\ge 0\;\;\;\forall\;x\in\Om,$$
then there is a point $q\in \Om$ such that
$$(u-v)(q)=\inf_{\Om}(u-v)\;\;\;\mbox{and}\;\;\;f(q,u(q))\le g(q,v(q)).$$
\end{lem}

\NI{\bf Proof.}  We employ the ideas in \cite{CIL} and use the concept of sub-jets and sup-jets.  We will prove part (a). The proof of part (b) will follow in an analogous manner. 
Set $M=\sup _{\Om}(u-v)$. By our hypothesis, $M>\sup_{\partial \Om}(u-v)$. Define, for $\vep>0$,
\be
w_{\vep}(x,y):=u(x)-v(y)-\frac{1}{2\vep}|x-y|^2,\;\;\;(x,y)\in \Om\times \Om.
\ee
Set $M_{\vep}:=\sup_{\Om\times \Om}w_{\vep}(x,y)$, and let $(x_{\vep},\;y_{\vep})\in \overline{\Om}\times\overline{\Om}$ be such that $M_{\vep}$ is attained at $(x_{\vep},\;y_{\vep}).$
The following are well-known, see \cite{CIL}.
\ben
 \lim_{\vep\rightarrow 0}M_{\vep}=\lim_{\vep\rightarrow 0}\left(u(x_{\vep})-v(y_{\vep})-\frac{|x_{\vep}-y_{\vep}|^2}{2\vep}\right)=M,\;\;\mbox{and}\;\;
 \lim_{\vep\rightarrow 0}\frac{|x_{\vep}-y_{\vep}|^2}{2\vep}=0.
 \een
Let $p\in \overline{\Om}$ be such that $x_{\vep}$ and $y_{\vep}\rightarrow p,$ as $\vep\rightarrow 0.$ Clearly, $M=u(p)-v(p).$ Since $M>\sup_{\partial \Om}(u-v)$, there is an  open set $O$, compactly contained in $\Om$, such that $p,\; x_{\vep}$ and $y_{\vep}\in O $.  

\NI Next, since $(x_{\vep},y_{\vep})$ is a point of maximum of $w_{\vep}(x,y)$, $((x_{\vep}-y_{\vep})/\vep, X_{\vep})\in \bar{J}^{2,+}u(x_{\vep})$ and $((x_{\vep}-y_{\vep})/\vep, Y_{\vep})\in \bar{J}^{2,-}v(y_{\vep})$. Moreover, we have, see \cite{CIL},
 \ben
- \frac{3}{\vep}\left(\begin{array}{cc}
 I &   0   \\
  0& I     \end{array}\right)
  \le
 \left(\begin{array}{lr}
 X_{\vep}&  0   \\
  0& -Y_{\vep} \end{array}\right)
  \le
\frac{3}{\vep} \left(\begin{array}{rr}
 I&   -I   \\
  -I& I
\end{array}\right).
\een
The above clearly implies $X_{\vep}\le Y_{\vep}$, and using the definitions of $\bar{J}^{2,+}$ and $\bar{J}^{2,-}$, we see that 
\eqRef{cin}
-f(x_{\vep},u(x_{\vep}))\le \left\langle \frac{X_{\vep}(x_{\vep}-y_{\vep})}{\vep},\frac{(x_{\vep}-y_{\vep})}{\vep}\right\rangle
\le\left\langle \frac{Y_{\vep}(x_{\vep}-y_{\vep})}{\vep},\frac{(x_{\vep}-y_{\vep})}{\vep}\right\rangle\le -g(y_{\vep}, v(y_{\vep})).
\end{equation}
Now let $\vep\rightarrow 0$ to conclude that $g(p,v(p))\le f(p,u(p)).\;\;\;\;\;\;\Box$
\vsp
We now state a few consequences of the above lemma. The first is an application of Lemma \ref{cmp0} to the eigenvalue problem (\ref{egpb}). This version will be used frequently in the rest of this work.

\begin{lem}\label{cmp}
Let $\lam_1$ and $\lam_2$ be real numbers, and $a(x)\in C(\Om)\cap L^{\infty}(\Om),\;a(x)>0$. Suppose that $u\in C(\overline{\Om})$ and $v\in C(\overline{\Om})$. \\
\NI (i) If $\sup_{\Om}(u-v)>\sup_{\p\Om}(u-v)$, and
$$\Df u+\lam_1 a(x) u^3\ge 0\;\;\;\;\mbox{and}\;\;\;\Df v+\lam_2 a(x) v^3\le 0,\;\mbox{in $\Om$},$$
then there is a point $p\in \Om$ such that $(u-v)(p)=\sup_{\Om}(u-v)$ and $\lam_1u^3(p)\ge \lam_2 v^3(p).$

\NI (ii) Similarly, if $\inf_{\Om}(u-v)<\inf_{\p\Om}(u-v)$, and 
$$\Df u+\lam_1 a(x) u^3\le 0\;\;\;\mbox{and}\;\;\;\Df v+\lam_2 a(x) v^3\ge 0,\;\;\mbox{in $\Om$,}$$ 
then there is a point $q\in\Om$ such that
$(u-v)(q)=\inf_{\Om}(u-v)$ and $\lam_1u^3(q)\le \lam_2 v^3(q).\;\;\Box $
\end{lem}
\vsp
We state below a consequence of Lemma \ref{cmp}. Versions of Lemma \ref{cmp1} are well-known in the context of eigenvalue problems for elliptic operators. Also see \cite{BNV,BD,J}. Here, we do not require that $\Om$ be bounded.

\begin{lem}\label{cmp1} Let $\Om\subset \IR^n$ be a domain. Suppose that $a(x)\in C(\Om),\;a(x)>0$, and $0<\lam_1<\lam_2$. Let $u\in C(\Om)$, and $v\in C(\Om), \;v>0,$ solve the problems
$$\Df u+\lam_1 a(x) u^3\ge 0\;\;\;\mbox{and}\;\;\;\Df v+\lam_2 a(x) v^3\le 0,\;\;\mbox{in $\Om$.} $$
Then either $u\le 0$ in $\Om$, or the following conclusions hold.\\
\NI (i) Let $U\subset \Om$ be a compactly contained sub-domain of $\Om$ such that $u>0$ somewhere in $U$. Then
$$\sup_{U} \frac{u}{v}=\sup_{\p U} \frac{u}{v}.$$
\NI (ii) Assume that $u>0$ somewhere in $\Om$. Suppose that $U_k\subset U_{k+1}\subset \Om,\;k=1,2,\cdots$, are compactly contained sub domains of $\Om$, with $\cup_k U_k=\Om$. If 
$\lim_{k\rightarrow \infty} \sup_{U_k}u/v<\infty,$ then 
$$\sup_{\Om}\frac{u}{v}=\lim_{k\rightarrow \infty} \left(\sup_{U_k}\frac{u}{v}\right)=\lim_{k\rightarrow \infty} \left(\sup_{\partial U_k}\frac{u}{v}\right).$$ 
\end{lem}
\NI{\bf Proof:} We prove (i). Let $U$ be a compactly contained sub-domain of $\Om$ and assume that $u>0$ somewhere in $U$. Suppose that $p\in U$ is such that $\sup_{U}(u/v)=u(p)/v(p)>\sup_{\p U} u/v.$ 
By our hypothesis, $u(p)>0$. Thus the function
\eqRef{cmp10}
w(x)=v(p)u(x)-u(p) v(x)\le 0,\;\;\;x\in \overline{U}.
\ee
In particular, $w(x)<0$ on $\p U$, and $w(p)=0$. Thus $\sup_U w>\sup_{\p U}w.$  Since $u(p)>0$ and $v(p)>0$, we have that for $\forall\;x\in \Om$,
$$\Df (v(p)u(x))+\lam_1a(x)(v(p)u(x))^3\ge 0,\;\mbox{and}\;\Df (u(p)v(x))+\lam_2a(x)(u(p)v(x))^3\le 0.$$
We may now apply Lemma \ref{cmp}(part(i)). It follows that there is a $z\in U$ such that $w(z)=\sup_U w$ and 
\eqRef{cmp11}
\lam_1a(z) u(z)^3 v(p)^3\ge \lam_2a(z) u(p)^3v(z)^3,
\ee
that is, $\tau u(p)/v(p) \le u(z)/v(z),$ where $\tau=\left(\lam_2/\lam_1\right)^{1/3}>1.$ This is a contradiction. 
Thus $\sup_{U}(u/v)=\sup_{\p U}(u/v).$

We now prove (ii).  Let $y\in \Om$ be such that $u(y)>0$. Take $k$ large, so that $y\in U_k$. Set $\mu_k=\sup_{\p U_k} (u/v)$. By part(i), the $\mu_k$'s are increasing. It is clear that the limit $\mu=\sup_k\mu_k<\infty.$ 
If $\sup_{\Om}(u/v)>\mu$ then one can find a set $U_k,$ for $k$ large, such that $\sup_{U_k}(u/v)>\mu$. This violates the maximum principle in part (i), as $\sup_{\p U_k}(u/v)\le \mu.$ 
The lemma holds.\quad $\Box$

\begin{rem}\label{cmp13}
As an application of Lemma \ref{cmp1}, we record the following. Let $\Om\subset \IR^n$ be a bounded domain, and $0<\lam_1<\lam_2$. Assume that 
$u,\;v\in C(\overline{\Om})$, $v>0,$ solve
$$\Df u+\lam_1 a(x) u^3\ge 0,\;\;\Df v+\lam_2 a(x) v^3\le 0,\;\;\mbox{for}\;x\in \Om,\;\;\mbox{and $v>u$ on $\p\Om$.}$$ 
Thus, if $u$ is positive somewhere in $\Om$ then $u$ is positive somewhere on $\p\Om.$
As a result, if $u\le 0$, on $\p\Om$, then $u\le 0$ in $\Om$.  $\Box$
\end{rem}

We now recall a few results from \cite{BMO1, BMO2, LW1, LW2} which we will utilize in our work. The first three lemmas contain versions of the comparison principle that apply in our context.

\begin{lem}\label{cmp2} 
Suppose that $f\in C(\Om)$, $f>0$, $f<0$ or $f\equiv 0$ in $\Om$. Let $u,v\in C(\overline{\Om})$ satisfy $\Df u\ge f(x)$ and $\Df v\le f(x)$ in $\Om$. Then
$$\sup_{\Om}(u-v)=\sup_{\p\Om}(u-v).\qquad \Box$$
\end{lem}

\begin{lem}\label{cmp3} 
Suppose that $f_1,\;f_2\in C(\Om)$ with $f_1(x)>f_2(x)$ in $\Om$. Let $u,v\in C(\overline{\Om})$ satisfy $\Df u\ge f_1(x)$ and $\Df v\le f_2(x)$ in $\Om$. Then
$$\sup_{\Om}(u-v)=\sup_{\p\Om}(u-v).\qquad \Box$$
\end{lem}

\begin{lem}\label{cmp30} 
Suppose that $f(x,t)\in C(\Om\times \IR,\IR)$ is strictly increasing in $t$. Let $u,v\in C(\overline{\Om})$ satisfy $\Df u\ge f(x,u)$ and $\Df v\le f(x,v)$ in $\Om$. If $u\le v$ on $\p\Om$ then $u\le v$ in $\Om$.\quad $\Box$
\end{lem}
The following estimate will prove useful in this work, see Theorem 5.1 in \cite{BMO2}. For a function $g$, define $g^+=\max\{g,0\}$ and $g^-=\min\{g,0\}$. Set $\s=3^{4/3}/4.$
\vsp
\begin{lem}\label{cmp4}
Let $\Om\subset \R^n$ be a bounded domain, and $B_{R_o}(z_o),\;z_0\in \IR^n,$ be the out-ball of $\Om$. 
Suppose $f\in C(\Om)\cap L^\infty(\Om)$, and $b\in C(\p\Om)$. If $u\in C(\overline{\Om})$ solves
$$\Df u=f(x),\;\;x\in\Om,\;\;\;u=b\;\;\;\mbox{on}\;\;\p\Om,$$ 
then the following bounds hold.
$$\inf_{\partial \Om}b-\s(\sup_{\Om} f^+)^{1/3}R_o^{4/3}\le u(x)\le   \sup_{\partial \Om}b-\s(\inf_{\Om} f^-)^{1/3}R_o^{4/3},\;\;x\in \Om.$$
 In particular, if $f(x)=-\lam a(x) u^3,\;a>0$, $\lam>0$, and $\mu=\sup_{\Om} a$, then a solution $u$ to (\ref{egpb}) satisfies
\ben
\inf_{\p\Om} b+\s (\lam \mu)^{1/3}R_o^{4/3}\inf_{\Om} u^- \le u(x)\le \sup_{\p\Om} b+\s(\lam \mu)^{1/3}R_o^{4/3} \sup_{\Om} u^+,\;\; x\in \Om.
\een
 Setting $\lam_0=(\s^3 \mu R_0^4)^{-1},$ then the above may be written more compactly as 
$$\inf_{\p\Om} b+(\lam/\lam_0)^{1/3}\inf_{\Om} u^- \le u(x)\le \sup_{\p\Om} b+\s(\lam/\lam_0)^{1/3}\sup_{\Om} u^+.\;\;\Box$$
\end{lem} 
\vsp
We also recall the following existence result proven in Theorem 3.1 in \cite{BMO2}, also see Corollary 3.3 and Theorem 5.5 therein. This will be used in showing the existence of solutions to equations related to the eigenvalue problem.

\begin{thm}\label{exst0}
Let $f\in C(\Om\times\R,\R)$ satisfy the condition $\sup_{\Om\times I}|f(x,t)|<\infty$, for any compact interval $I$, and $b\in C(\p\Om)$. Consider the following Dirichlet problem
$$(\star)\qquad \Df u=f(x,u(x)),\;\mbox{in $\Om$, and $u=b$ on $\p\Om$.}$$
(a) Suppose that \\
\NI (i) $u_*\in C(\overline{\Om})$ is a sub-solution of $(\star)$, i.e., $\Df u_*\ge f(x,u_*),$ in $\Om$, and $u_*\le b$ on $\p\Om$, and \\
\NI (ii) $u^*\in C(\overline{\Om})$ is a super-solution of $(\star)$, i.e., $\Df u^*\le f(x,u^*)$, in $\Om$,
and $u^*\ge b$ on $\p\Om$.\\
 If $u_*\leq u^*$ in $\Om$ then problem $(\star)$ admits a solution $u\in C(\overline{\Om})$ such that $u_*\leq u\leq u^* $ in $\Om$. \\
\NI (b) If $f$ is such that any solution to $(\star)$ has a priori supremum bounds, then there is a solution $u\in C(\overline{\Om})$ to $(\star)$.
\quad$\Box$
\end{thm}
\vsp
We now record a local Lipschitz continuity result, proven in \cite{BMO2}, see Theorem 2.4 therein. Also see \cite{LW2}.
\begin{lem}\label{lpc}
Let $\al$ be a constant. Any solution $u\in C(\Om)\cap L^{\infty}(\Om)$ of $\Df u(x)  \geq  \al,$ in $\Om$, is locally Lipschitz continuous in $\Om$. More specifically, given $x_0\in\Om$ there is a constant $C$ that depends on $x_0,\,\mbox{diam}(\Om),\;|\al|$ and $\|u\|_{L^\infty(\Om)}$ such that
$$|u(x)-u(y)|\leq C|x-y|,\;\;\;\;\;x,y\in B_{d}(x_0),$$
where $d:=\mbox{dist}(x_0,\p\Om)/3$. A similar result holds if $\Df u\le \al$ in $\Om$. \quad$\Box$
\end{lem}

We now shift our attention to obtaining estimates for a problem that is related to (\ref{egpb}). These will be important in proving the existence of the first eigenvalue and an associated eigenfunction. 
To achieve this purpose, we study the following Dirichlet problem. Let $a\in C(\Om)\cap L^{\infty}(\Om),\;a>0$, $\dl>0$ and $\lam>0$. Consider positive solutions to 
the problem
\eqRef{Frd}
\Df u+\lam a(x) u^3(x)=0\;\mbox{in $\Om$, and $u=\dl$ on $\p\Om$}.
\end{equation}
In order to show existence we note that the function $\psi=\dl$ is a sub-solution to (\ref{Frd}). For small $\lam$, we obtain a priori supremum bounds. This will lead to
the existence of a solution $u$. 

\begin{lem}\label{ap} Let $\Om\subset \IR^n$ be a bounded domain. Suppose that $a(x)\in C(\Om)\cap L^{\infty}(\Om),\;a(x)>0$, $\dl\ge 0$, and $\lam>0$. Let $R_o$ be the radius of 
the out-ball for $\Om$, $\mu=\sup_{\Om} a$, $\s=3^{4/3}/4$ and $\lam_0=(\s^3 \mu R_o^4)^{-1}.$ Consider the problem
$$\qquad (\star)\quad\Df u+\lam a(x) u^3(x)=0,\;\mbox{in $\Om$, and $u=\dl$ on $\p\Om$}.$$
Assume that $u\in C(\overline{\Om})$ is a solution to $(\star)$.\\
 \NI (i) If $\lam=0$ then $u=\dl$ in $\Om$.\\
\NI (ii) If $\lam<0$ and $\dl>0$ then $0\le u<\dl$. If $\dl=0$ then $u=0$ is the only solution. \\
\NI (iii) If $\dl=0$ and $u\in C(\overline{\Om})$ is a non-trivial and non-constant solution, then $\lam> 0$. \\
\NI (iv) If $0<\lam<\lam_0$ then $u$ is positive in $\Om$ and a priori bounded. More precisely,
$$\dl<u\le \sup_{\Om}u\le \frac{\dl}{1-(\lam/\lam_0)^{1/3}}.$$
\end{lem}

\NI{\bf Proof:} We show (i). If $\lam=0$ then $u$ is infinity-harmonic and $u=\dl$ in $\Om$. For part (ii), suppose that $\lam<0$. Let $\Om^-$ denote the set where $u<0$. Then $\Df u=|\lam| a(x) u^3\le 0$, in $\Om^-$, with $u$ vanishing on $\p\Om^-$. But $u>0$, in $\Om^-$, since $u$ is infinity super-harmonic in $\Om^-$. It follows that
$\Om^-=\emptyset$ and $u\ge 0$, in $\Om$. Thus, $u$ is infinity sub-harmonic in $\Om$, and $0\le u<\dl$. If $\dl=0$, we get $u=0$ in $\Om$, for $\lam\le 0.$  Clearly, parts (i) and (ii) imply part (iii).

We now prove part (iv). We will assume that $\dl>0$, the conclusion for $\dl=0$ follows quite easily. We recall Lemma \ref{cmp4},
\eqRef{est0}
\dl+\left(\frac{\lam}{\lam_0} \right)^{1/3} \inf_{\Om} u^-\le u(x)\le \dl +\left(\frac{\lam}{\lam_0}\right)^{1/3}\sup_{\Om} u^+.
\end{equation}
If $\inf_{\Om} u^-<0$, then (\ref{est0}) leads to
$$\frac{\dl}{1-(\lam/\lam_0)^{1/3}}\le \inf_{\Om} u^-<0,$$
a contradiction. Thus (\ref{est0}) yields
$$0\le u\le \sup_{\Om}u\le \frac{\dl}{1-(\lam/\lam_0)^{1/3}}.$$
Since $u$ is infinity super-harmonic, $u>\dl$ in $\Om$. \quad $\Box$
\vsp
Finally, we prove that nontrivial solutions to (\ref{Frd}), when $\dl\ge 0$, change sign for large enough $\lam.$ This was first shown in \cite{BMO2} and implies that, in the event eigenfunctions corresponding to large eigenvalues exist, these eigenfunctions would change sign, a fact well-known for the case of elliptic operators. Its relevance to our current work is in obtaining lower and upper bounds for the first eigenvalue. We provide a proof of this result for completeness. We do not assume that $a(x)>0$ everywhere in $\overline{\Om}$. 

\begin{thm}\label{eqn1}
Let $\Om\su \R^N$ be a bounded domain, and
$a(x)\in C(\Om)\cap L^{\infty}(\Om),\;a(x)\ge 0,$ and
$a(x)\not\equiv 0.$ Set $\mu=\sup_{\Om} a$, $\s=3^{4/3}/4$ and $\lam_0=(\s^3 \mu R_o^{4})^{-1}$, where $R_o$ is the radius of the out-ball for $\Om$.
For $0<\al<1$, define $\Om_{\al}=\{x\in\Om:\;a(x)>\al \mu\}$, and set $\rh_{\al}$ to be the radius of the in-ball of $\Om_{\al}.$ 
Let $\dl\ge 0$, 
and suppose that $(\lam,u),\;u\not\equiv0$, solves 
\eqRef{eqn2}
\Df u+\lam a(x) u^3=0,\;x\in \Om,\;\;\mbox{and}\;u=\dl\;\;\mbox{on $\p\Om$}.
\ee
Set $$\Lambda=  \frac{4^4}{3^3\s^3\mu}\left(\inf_{0<\al\le 1}\left(\frac{1}{\al\rho_{\al}^4}\right)\right)<\infty.$$
(i) If $\dl=0$, then $\lam\ge \lam_0$.\\
(ii) If $\dl> 0$ and $u\ge 0$, then we have the upper bound $\lam<\Lambda.$ If $\dl=0$ and $u\ge 0$ then $\lam_0\le \lam\le \Lambda.$
In any case, if $\lam$ is large enough then every solution $u$ to (\ref{eqn2}) changes sign in $\Om$, regardless of $\dl.$
 \end{thm}

\noindent {\bf Proof.} For part (i), we refer to parts (iii) and (iv) of Lemma \ref{ap}. By (\ref{est0}), if $\dl=0$ then $u=0$, for $\lam<\lam_0.$ 

We now prove part (ii). See Lemma \ref{ap} (ii) and Theorem \ref{exst0} (b) for the lack of a lower bound for $\lam$ when $\dl>0$. If $\dl=0$ and $\lam\ge \lamo$ then $u>0$, since $u$ is infinity super-harmonic. 

In order to show the upper bound for $\lam$, we assume that $\lam>0$. Let $(\lam,u)$, $u\in C(\overline{\Om}),\;u>0$, solve (\ref{eqn2}).
Being infinity super-harmonic in $\Om$, $u$ satisfies the strong minimum principle  and $u>\dl$.
For $0<\al<1$, let $B_{\rho_{\al}}(z_{\al})$ be the in-ball for $\Om_{\al}$. For $0\le r\le \rho_{\al}$, define $m(r)=\inf _{\p B_r(z_{\al})}u$. Then $\dl<m(r)\le u$, in $B_r(z_\al)$, and $m(r)$ is decreasing.
Consider 
$$v(x)=\dl+(m(0)-\dl)\left(1-\frac{|x-z_{\al}|}{\rho_{\al}}\right),\;\;x\in B_{\rh_{\al}}(z_{\al}).$$
It is clear that $v$ is infinity harmonic in $B_{\rho_{\al}}(z_{\al})\setminus \{z_{\al}\}$. Since $u\ge \dl$ on $\p B_{\rho_{\al}}(z_{\al})$ and $u(z_{\al})=m(0)$, by Lemma \ref{cmp2}, 
$v\le u$ in $B_{\rho}(z_{\al})\setminus\{z_{\al}\}$. Taking $|x-z_{\al}|=\tht \rh_{\al}$, for $0\le \theta <1$, and noting that $v(\tht\rh_{\al})\le m(\tht\rh_{\al})$, we have
\begin{equation}\label{mro} 
\frac{m(0)-\dl}{m(\theta\rho_{\al})-\dl}\le \frac{1}{1-\theta}.    
\end{equation}
Next we consider, in the ball  $B_{\theta \rho_{\al}} (z_{\al})$, the function
$$
w(x)=\s(\al \mu\lam)^{1/3}m(\theta\rho_{\al})\left( \left(\theta\rho_{\al}\right)^{4/3}-|x-z_{\al}|^{4/3}\right)+m(\theta\rho_{\al}).
$$
Using (\ref{radl}), a calculation shows that
$$
\Delta_{\infty}w=-\al\lam \mu m(\theta\rho_{\al})^3,\; \mbox{in $B_{\theta\rho_{\al}}(z_{\al}),$ and $  w=m(\theta\rho_{\al})$ on $|x-z_{\al}|=\theta\rho_{\al}$}.
$$
In $B_{\theta\rho_{\al}}(z_{\al})\subset \Om_{\al}$, we note that $a(x)>\al \mu$ and $u>m(\theta \rho_{\al})$. Thus
$$\Delta_{\infty}u=-\lam a(x) u^3<-\al\lam \mu m(\theta\rho_{\al})^3, \;\;\;x\in B_{\theta\rho_{\al}}(z_{\al}),$$
with $u\ge w$ on $|x-z_{\al}|=\theta\rho_{\al}$. Lemma \ref{cmp3} yields that
$w\le u$ in $B_{\theta\rho_{\al}}(z_{\al})$. Moreover,
$$w(z_{\al})=\s\lam^{1/3}(\al \mu)^{1/3}m(\theta\rho_{\al})\left(\theta\rho_{\al}\right)^{4/3}+m(\theta\rho_{\al})\le u(z_{\al})=m(0).$$
Recalling that $u>\dl$ and rewriting,
$$\s\lam^{1/3}(\al \mu)^{1/3}(m(\theta\rho_{\al})-\dl)\left(\theta\rho_{\al}\right)^{4/3}+m(\theta\rho_{\al})-\dl\le m(0)-\dl.$$
Rearranging and using (\ref{mro}), we have 
$$ \s\lam^{1/3}(\al \mu)^{1/3}\left(\theta\rho_\al\right)^{4/3}+1\le \frac{m(0)-\dl}{m(\theta\rho_{\al})-\dl}\le \frac{1}{1-\theta}.$$   
Rewriting, we get
$$\s \lam^{1/3}(\al \mu)^{1/3}\rho_\al^{4/3}\le \frac{1}{\theta^{1/3}(1-\theta)},\;\;\;0<\theta<1.$$
By computing the minimum of the right hand side, which occurs at $\theta=1/4$, we obtain
$$
\lam\le\frac{4^4}{3^3\s^3\mu(\al \rho_{\al}^4)}.
\eqno{\qed} $$

\section{\bf Existence and properties of solutions to (\ref{Frd})}

In this section, we derive properties of solutions to (\ref{Frd}) when $u$ takes positive values on $\p\Om.$ This will lead to an existence result for (\ref{Frd}) with non-trivial right-hand side. All these will be proven under the condition that $\lam$ is less than the first eigenvalue $\lamo$ of $\Df $.  We will adapt the comparison principle in Lemma \ref{cmp1} to the current context and this will lead to uniqueness, under some conditions.  

We will begin with a discussion of how to define the first eigenvalue. The  basic idea resembles closely the one employed in \cite{BNV, BD, J}. 

\begin{lem}\label{exst2} Let $\Om\subset \IR^n$ be a bounded domain. Suppose that $a(x)\in C(\Om)\cap L^{\infty}(\Om),\;a(x)>0,$ and assume that $\dl>0$. Define $\lam_0=(\s^3 \mu R_o^4)^{-1}$, where
$\s=3^{4/3}/4$, $\mu=\sup_{\Om} a$, and $R_o$ the radius of the out-ball of $\Om$.
Then the Dirichlet problem
\eqRef{pbst}
\Df u+\lam a(x) u^3=0,\;\mbox{in $\Om$, and $u=\dl$ on $\p\Om$},
\ee
has a positive solution $u$ for $0\le \lam<\lam_0$.
\end{lem} 
\NI{\bf Proof.} We use  Theorem \ref{exst0}(b) and Lemma \ref{ap}(iv). Since $\lam<\lam_0$, any solution $u$ is a priori bounded and 
Theorem \ref{exst0} leads to a solution. Lemma \ref{ap} ensures that $u>\dl$ in $\Om$. \quad $\Box$
\vsp
We now discuss the definition of the first eigenvalue. The fact that it is indeed an eigenvalue and has at least one eigenfunction will be shown in Section 4. We define, for each $\dl>0$,
\eqRef{set}
S=S(\Om)=\{\lam\ge 0:\;\mbox{Problem (\ref{Frd})(or (\ref{pbst})) has positive solutions}\}.
\end{equation}
By Lemma \ref{exst2}, $S$ is non-empty. By Theorem \ref{eqn1}, $S$ is bounded above. Now set
\eqRef{set1}
\lamo=\sup_S \lam.
\ee 
We refer to $\lamo$ as the first or the principal eigenvalue of $\Df$ on $\Om$. 

\begin{rem}\label{set21} We record the following conclusions.\\
\NI{\bf (i)} By Lemma \ref{ap}, $\lamo\ge (\s^3\mu R_0^4)^{-1}.$ We show that the interval $[0,\lamo) \subset S$. Let $\lam\in S$ and $u>0$ be a solution to
$$\Df u+\lam a(x) u^3=0,\;\mbox{in $\Om$, and $u=\dl$ on $\p\Om$}.$$
Note that $u>\dl$ in $\Om$. If $0<\lam'<\lam$, then $u$ is a super-solution to
\eqRef{lft}
\Df v+\lam' a(x) v^3= 0,\;\mbox{in $\Om$, and $v=\dl$ on $\p\Om$}.
\end{equation}
Clearly, $w=\dl$ is a sub-solution; it follows from Theorem \ref{exst0} that there is a solution $v$ to (\ref{lft}) such that $\dl<v\le u.$ Hence, $\lam' \in S$. That $\lamo\not\in S$ will follow from Lemma \ref{up} below. \\
\NI{\bf (ii)} The set $S$ is independent of the value of $\dl.$ This follows by scaling.\\
\NI{\bf (iii)} We discuss the influence of the weight function $a(x)$. Write in (\ref{set}), $S=S(\Om,a)$ and in (\ref{set1}), $\lamo=\lamo(a)$.  We claim that $S(\Om,b)\subset S(\Om,a),$ and $\lamo(a)\ge \lamo(b),$ when $0\le a(x)\le b(x)$, in $\Om$. 
To see this, let $\lam\in S(\Om,b)$. We can find a function $u\in C(\overline{\Om}),\;u>0,$ that solves $\Df u+\lam b(x) u^3=0,$ in $\Om$, and $u=\dl$ on $\p\Om$. Then 
$\Df u+\lam a(x) u^3\le 0,$ in $\Om$. Since $v=\dl$ is a sub-solution, we have from Theorem \ref{exst0} that there is a function $\bar{u}\in C(\overline{\Om}),$ $v\le \bar{u}\le u$, that solves 
$$\mbox{$\Df \bar{u}+\lam a(x)\bar{u}^3=0$, in $\Om$, and $\bar{u}=\dl$ on $\p\Om$.}$$
Thus $\lam\in S(\Om,a)$ and $\lam_{\Om}(a)\ge \lamo(b).$\\
\NI{\bf (iv)} By Theorem \ref{eqn1}, the set $S$ is bounded from above and $\lamo<\infty$. \quad $\Box$
\end{rem}

Later in this section, we will use (\ref{set21}) to state an existence result for boundary data that has one sign, under the hypothesis $0\le \lam<\lamo$. A related result is in Lemma \ref{nst} where it is shown that if $0\le \lam<\lamo$ and the boundary data is zero then the zero solution is the only solution. 

We restate problem (\ref{Frd}) for easy reference. Also recall (\ref{set}) and (\ref{set1}).
We  will study the properties of a solution $u\in C(\overline{\Om}),\;u>0,$ to 
\eqRef{Frd0}
\Df u+\lam a(x) u^3=0,\;\mbox{in $\Om$, and $u=\dl>0$ on $\p\Om$.}
\end{equation}
Here $0<\lam\le \lamo<\infty.$  We refer the reader to Lemma \ref{ap} for the case $\lam\le 0$.

We show next that if $\lam\in S,$ then, for some $\vep>0$, $\lam+\vep$ is also in $S$. This will imply that $\lamo\not\in S$, justifying part  (iii) in Remark \ref{set21}. 
\vsp
\begin{lem}\label{up} Let $a(x)\in C(\Om)\cap L^{\infty}(\Om)$ with $a(x)>0$. Suppose that for some $\lam>0$, there is a function $v\in C(\overline{\Om}),\;v>0,$ such that 
\eqRef{up1}
\Df v+\lam a(x) v^3\le 0,\;\mbox{in $\Om$, and $v\ge \dl$ on $\p\Om$.}
\end{equation}
 Set $m=\sup_{\Om} v$. Then, for every $\vep$ such that $0<\vep<\lam(\dl/m)^3$ the problem 
$$\Df u+(\lam+\vep)a(x) u^3=0,\;\mbox{in $\Om$, and $u=\dl$ on $\p\Om$,}$$
has a positive solution $u\in C(\overline{\Om})$. Hence, $\lamo\not\in S$, where $S$ is as in  (\ref{set}).
\end{lem}

\NI{\bf Proof:} We apply Theorem \ref{exst0} to achieve the proof.
Let $0<\vep<\lam(\dl/m)^3$. Take $0<\al<1$ such that 
\eqRef{up0}
\vep<\al\lam\left(\dl/m\right)^3.
\ee
Since $v>0$, it follows that $v>\dl$ in $\Om$. Define
$$w(x)=v(x)-\al \dl,\;\;\;x\in \Om.$$
Then (\ref{up1}) becomes
\eqRef{up2}
\Df w+\lam a(x)v^3\le 0,\;\mbox{in $\Om$,}\;\;\mbox{and $w\ge (1-\al)\dl$ on $\p\Om$.}
\end{equation}
Writing $v=w+\al\dl$ and noting $w\le m$ in $\Om$, we expand, using (\ref{up2}), to obtain 
\bea 
\Df w+(\lam+\vep)a(x)w^3&\le& a(x)((\lam+\vep)w^3 -\lam v^3)\nonumber  \\
&=&a(x)\left((\lam+\vep) w^3-\lam (w+\al\dl)^3\right)\nonumber\\
&\le& a(x)\left(\vep m^3-\lam(3\al \dl w^2+3\al^2\dl^2w+\al^3\dl^3)\right).\label{up4} 
\eea  
Since $w\ge (1-\al)\dl$ and $\al^2-3\al+3>1$, for $0<\al<1$, we have that 
$$
\al\dl\left( 3w^2+3\al\dl w+\al^2\dl^2\right)\ge \al\dl^3(\al^2-3\al+3)>\al\dl^3.
$$
Using the above in (\ref{up4}) and applying (\ref{up0}),
\eqRef{up5}
\Df w+(\lam+\vep)a(x) w^3\le a(x)\left(\vep m^3-\lam\al\dl^3\right)\le 0,\;\;x\in \Om.
\end{equation}
 It is clear that if we take $0<\vep<\lam(\dl/m)^3$ and any $\al$ with $(\vep/\lam)(m/\dl)^3<\al<1$ (see (\ref{up0})) then the function
\eqRef{sups}
h=h(\al)=\frac{w}{(1-\al)}=\frac{v-\al \dl}{1-\al}\ge \dl,
\end{equation} 
 defined in $\Om$, is a super-solution to
\eqRef{up6}
\Df f+(\lam+\vep)a(x) f^3=0,\;\mbox{in $\Om$, and $f=\dl$ on $\p\Om.$}
\ee
Next, we observe that the function $g(x)=\dl,\;x\in \overline{\Om}$ is a sub-solution of (\ref{up6}). Since $g\le h$ in $\Om$, invoking Theorem \ref{exst0}, we obtain that 
(\ref{up6}) has a solution $u$ such that $g\le u\le h$ in $\Om$.
\quad $\Box$

We prove now a comparison principle by employing Lemmas \ref{cmp1} and \ref{up}. This will imply the uniqueness of solutions to (\ref{Frd0}) for $0\le \lam<\lamo$. We will utilize the function $h$ defined in (\ref{sups}). Also, see \cite{BNV, BD, J}. We do not assume that $\Om$ is bounded. 

\begin{lem}\label{ord3}
Suppose that $a(x)\in C(\Om)\cap L^{\infty}(\Om),\;a(x)>0$ and $\lam>0$. Let $u,\;v\in C(\Om),\;v>0,$ solve the problems
$$\Df u+\lam a(x) u^3\ge 0\;\;\mbox{and}\;\;\Df v+\lam a(x) v^3\le 0\;\;\mbox{in $\Om$}.$$
Either $u\le 0$ in $\Om$, or the following holds. \\
\NI (a) If $U$ is a compactly contained sub-domain of $\Om$ and $u>0$ somewhere in $U$, then $\sup_U(u/v)=\sup_{\partial U}(u/v).$ \\
\NI (b) Suppose that $u>0$ somewhere in $\Om$ and $\{U_m\},\;m=1,2,\cdots$ is an increasing sequence of compactly contained sub-domains of $\Om$, with $\cup_{m=1}^{\infty}U_m=\Om$. If 
$\lim_{m\rightarrow \infty}\sup_{\partial U_m}(u/v)=k<\infty$, then $k>0$ and $u\le k v$ in $\Om$.
\end{lem}

\NI{\bf Proof:} We take part (a). Let $U$, compactly contained in $\Om$, be such that $u>0$ somewhere in $U$. Set $\ell=\ell(U)=\inf_{\p U} v$. Being infinity super-harmonic, $v>\ell$ in $U$. 
If we define, for $0<\al<1$, 
$$h=\frac{v-\al \ell}{1-\al},\;\;\mbox{in $U$,}$$ 
then a simple calculation shows that $h\ge v\ge \ell(1-\al)$ in $\overline{U}$. By (\ref{up5}), we also have
$$\Df h+(\lam+\vep)a(x) h^3\le 0,\;\mbox{in $U$},$$ 
where $0<\vep<\lam \al (\ell/\sup_{U} v)^3.$ By Lemma \ref{cmp1}, we have that for every $0<\al<1$,
$$\sup_{U}\frac{u}{h}= \sup_{\p U} \frac{u}{h}.$$
Letting $\al\downarrow 0$, we obtain that 
$$\sup_U \frac{u}{v}=\sup_{\p U}\frac{u}{v}.$$
Part (b) of the lemma follows by applying the arguments of Lemma \ref{cmp1} and Remark \ref{cmp13}.\quad $\Box$
\vsp
As a consequence of Lemma \ref{ord3}, we obtain the uniqueness of solutions to the Dirichlet problem (\ref{Frd0}) with $\inf_{\p\Om} b>0.$ 

\begin{rem}\label{unq} Let $\Om\subset \IR^n$ be a domain, $0<\lam<\infty$, and $a(x)\in C(\Om)\cap L^{\infty}(\Om),\;a(x)>0$. Suppose that $u,\;v\in C(\overline{\Om})$, with $u>0$ and $v>0$, solve
\ben
\Df u+\lam a(x) u^3=0\;\;\;\mbox{in $\Om$, and}\;\;\;\;\Df v+\lam a(x) v^3=0.
\een
Suppose that $\{U_m\},\;m=1,2,\cdots,$ is an increasing sequence of compactly contained sub-domains of $\Om$ with $\cup_{m=1}^{\infty} U_m=\Om$.
If $\lim_{m\rightarrow \infty} (\sup_{\p U_m}(u/v))\;\mbox{and}\;\lim_{m\rightarrow \infty} (\sup_{\p U_m}(v/u))$
exist then
$$\lim_{m\rightarrow \infty} \left(\inf_{\p U_m}\frac{u}{v}\right)\le \frac{u(x)}{v(x)}\le \lim_{m\rightarrow \infty} \left(\sup_{\p U_m}\frac{u}{v}\right),\;\;\;x\in \Om.$$
These limits, if they exist, are independent of the sequence.

As an application, if $\Om$ is bounded, $u,\;v\in C(\overline{\Om})$, $b\in C(\p \Om)$ is such that $\inf_{\p\Om} b>0$ and $u=v=b$ on $\p\Om$, then we have that $u=v$ in $\Om$.  
\quad $\Box$
\end{rem}
\vsp
Next we record an application of Lemma \ref{cmp1}. This will be used in Section 4, where we show the existence of the first eigenvalue.

\begin{rem}\label{set3} Let $0<\lam<\lam^{\prime}$. Suppose that $(\lam,u),\;u>0,$ and $(\lam^{\prime},v),\;v>0,$ solve the problem (\ref{Frd0}). As $u$ and $v$ take the same boundary data, by Lemma \ref{cmp1}, $u\le v$ in $\Om$. Thus, if $\lam_k\uparrow \lamo$ then the corresponding unique solutions $ \{v_k\}$ form an increasing sequence. $\Box$
\end{rem}

We now show that if $\dl=0$ in (\ref{Frd0}) and $\lam<\lamo$, then the  only solution is the zero solution.
 The proof requires the existence of a solution that is positive in $\overline{\Om}$. Note that this is guaranteed by the nature of the set $S$, see (\ref{set21}). 
\vsp
\begin{lem}\label{nst} Let $a(x)\in C(\Om)\cap L^{\infty}(\Om),\;a(x)>0$, and $\lam>0$. Suppose that $v\in C(\overline{\Om}),\;v>0,$ and $v$ solves
$$\Df v+\lam a(x) v^3\le 0,\;\;\mbox{in $\Om$.}$$
If $\inf_{\p\Om} v>0$, and $u\in C(\overline{\Om})$ solves
\eqRef{nst10}
\Df u+\lam a(x) u^3\ge 0,\;\mbox{in $\Om$, and $u=0$ on $\p\Om$},
\ee
then $u\le 0$ in $\Om$. If equality holds in (\ref{nst10}) then $u=0$ in $\Om.$
\end{lem}

\NI{\bf Proof:} We use Lemma \ref{ord3}. If $u$ solves (\ref{nst10}) and $u$ is positive somewhere in $\Om$ then
$\sup_{\Om}(u/v)= \sup_{\p\Om}(u/v)>0.$ This being a contradiction,  we have $u\le 0$ in $\Om$.
If, instead of the inequality in (\ref{nst10}), equality holds,  then both $u$ and $-u$ are solutions. We conclude that $u=0$ in $\Om$.
Incidentally, if $0<\lam<\lamo$ then such a function $v$ exists, by (\ref{set21}). 
\quad $\Box$
\vsp
A related result follows below.

\begin{rem}\label{nst1} Let $a(x)\in C(\Om)\cap L^{\infty}(\Om),\;a(x)>0$, $0<\lam_1<\lam_2$ and $\dl>0$. Suppose that $u,\;v\in C(\overline{\Om})$ solve the problems
$$\Df u+\lam_1 a(x) u^3\ge 0\;\;\mbox{and}\;\;\Df v+\lam_2 a(x) v^3\le 0,\;\mbox{in $\Om$.}$$
Assume also that $u\le \dl\le v$ on $\p\Om$.
Set $\al=(\lam_1/\lam_2)^{1/3}$. If $v>0$ and $u$ is positive somewhere in $\Om$, then we claim that
$$u(x)\le \dl^{1-\al}(sup_{\Om} v^{\al}),\;\;\;\forall\;x\in \Om.$$ 
To see this, we make the following observation. Let $\lam>0$ and $w\in C(\overline{\Om})$ be positive. If $\Df w+\lam a(x) w^3\ge 0$, in $\Om$, then for any $\bt>1$, we have that
$\Df w^{\bt}+\lam \bt^3a(x) w^{3\bt}\ge0$. If instead, $\Df w+\lam a(x) w^3\le 0$, then for any $0<\bt<1$, it follows that $\Df w^{\bt}+\lam \bt^3a(x) w^{3\bt}\le 0$.

Now take $\beta=\al$. Since $\al<1$, we invoke Lemma \ref{ord3} to conclude that 
$u/v^{\al}\le  \sup_{\p\Om}(u/v^{\al}) =\dl^{1-\al}.$
The claim holds. 

We surmise that a stronger estimate holds, namely, that $u(x)\le C\dl,\;\forall\;x\in \Om$, where $C=C(\lam_1,\lam_2,\Om)$. However, a proof is not yet clear to us.
$\quad\Box$
\end{rem}

We now state the first of the two existence results of this section. We include a partial result about uniqueness. Also see \cite{J}.

\begin{thm}\label{exst3}
Let $\Om\subset \IR^n$ be a bounded domain, and $a(x)\in C(\Om)\cap L^{\infty}(\Om)$ with $a(x)\ge 0$. Suppose that $0\le \lam<\lamo$ and $b\in C(\p\Om)$. 
Then there is a function $u\in C(\overline{\Om})$ that solves the following Dirichlet problem, that is,
\eqRef{exst4}
\Df u+\lam a(x) u^3=0,\;\mbox{in $\Om$, and $u=b$ on $\p\Om$.}
\ee
In addition, we have the following. \\
\NI (i) If $b= 0$ on $\p\Om$, then $u=0,$ in $\Om$. \\
\NI (ii) Suppose that $b\not= 0$ on $\p\Om$. If $\inf_{\p\Om}b\ge 0$ or $\sup_{\p\Om}b\le 0$ then every solution 
$u\in C(\overline{\Om})$ is non-vanishing in $\Om$. \\
\NI (iii) If $\inf_{\p\Om} |b|>0$ then $u$ is unique.
\end{thm}

\NI{\bf Proof.} We first show the existence of a solution to (\ref{exst4}). Let $m=\sup_{\p\Om}b$ and $\ell=\inf_{\p\Om} b.$ If $\ell=m$, Remark \ref{set21} 
gives us a solution. Take
$m_1>\max(m,0)$, and $\ell_1<\min(0,\ell).$ \\
\NI By Remark \ref{set21}, there is a $w_1\in C(\overline{\Om}),\;w_1>0,$ that solves
\eqRef{exst5}
\Df w_1+\lam a(x) w_1^3=0,\;\mbox{in $\Om$, and $w_1=m_1$ on $\p\Om$.}
\ee
By (\ref{exst5}), the function $w_2=(\ell_1/m_1)w_1$ solves
$$\Df w_2+\lam a(x) w_2^3=0,\;\mbox{in $\Om$, and $w_2=\ell_1$ on $\p\Om$.}$$
Clearly, $w_2\le w_1$, in $\Om$, and $w_2\le b\le w_1$ on $\p\Om$. 
By Theorem \ref{exst0}, there is a solution $u\in C(\overline{\Om})$ to (\ref{exst4}) such that
$w_2\le u\le w_1$. 

It is clear that part (i) of the lemma follows from Remark \ref{set21} and Lemma \ref{nst}. We prove part (ii).
We will assume that $b\ge 0$  (if $b\le 0$, we work with $-u$). Suppose that $u$ changes sign in $\Om$.  Call $\Om^-=\{u<0\}$. Then $u$ solves
$$\Df u+\lam a(x) u^3=0,\;\mbox{and $u=0$ on $\p\Om^-$.}$$
Since $\lam<\lamo$, by Remark \ref{set21},  there is a solution $v\in C(\overline{\Om})$, for $\dl>0$, to 
$$\Df v+\lam a(x) v^3=0,\;v>0,\;\mbox{in $\Om$, and $v=\dl$ on $\p\Om$.}$$
Since $v\ge \dl,$ in $\overline{\Om^-}$, applying Lemma \ref{nst} to $u$ and $v$ in $\Om^-$, we obtain a contradiction.  Thus, $u\ge 0$ in $\Om$, and being 
infinity super-harmonic
we have that $u>0$ in $\Om$. Part (iii) follows from Remark \ref{unq}, also see Lemma \ref{ord3}. 
$\Box$
\vsp
We now state an existence result for non-homogenous right hand sides. We will prove this under the somewhat restrictive assumption
that $\inf _{\Om}a(x)>0$. We do not address the issue of uniqueness. We borrow an idea from Lemma \ref{up}. Also see \cite{J}.
\vsp
\begin{thm}\label{exst6}
Let $\Om\subset \IR^n$ be a bounded domain, $a(x)\in C(\Om)\cap L^{\infty}(\Om)$, with $\inf_{x\in\Om}a(x)> 0$, and $0\le \lam<\lamo$. Suppose that $h\in C(\Om)\cap L^{\infty}(\Om)$ and $b\in C(\p\Om)$. 
Then there is a function $u\in C(\overline{\Om})$ that solves the following Dirichlet problem, 
\eqRef{exst7}
\Df u+\lam a(x) u^3=h(x),\;\mbox{in $\Om,$ and $u=b$ on $\p\Om$.}
\ee
\end{thm}
\NI {\bf Proof.}  Our approach is similar to Lemma \ref{exst3}. Let $m=\sup_{\p\Om} b,\;\ell=\inf_{\p\Om} b$, $M=\sup_{\Om}|h|$ and $\nu=\inf_{\Om} a$. 
Take $m_1>\max(m,0)$ and $\ell_1<min(0,\ell, -m_1)$.  We will construct a sub-solution and a super-solution to (\ref{exst7}). 

(i) We first construct a super-solution. Let $w_1\in C(\overline{\Om}),\;w_1>0,$ be a solution to 
$$\Df w_1+\lam a(x) w_1^3=0,\;\mbox{in $\Om$, with $w_1=m_1$ on $\p\Om$.} $$
Existence follows from Remark \ref{set21}. Being infinity super-harmonic, $w_1>m_1$. For $0<\al<1$, take $w_2=w_1-\al m_1$. 
Thus $\Df w_2+\lam a(x) (w_2+\al m_1)^3=0.$ Expanding, 
$$\Df w_2+\lam a(x) w_2^3=-\lam a(x)\left( 3\al m_1 w^2_2+3\al^2m_1^2 w_2+\al^3m_1^3\right).$$
Noting that $w_2\ge (1-\al)m_1$, in $\Om$, we obtain that
$$\Df w_2+\lam a(x) w_2^3\le -\lam \nu m_1^3\left( 3\al (1-\al)^2+3\al^2(1-\al) +\al^3\right).$$
Set $w=w_2/(1-\al)$. Selecting $\al$ close enough to $1$, we obtain from above  that
$$\Df w+\lam a(x) w^3\le -\lam \nu m_1^3\left( \frac{3\al}{(1-\al)}+\frac{3\al^2}{(1-\al)^2} +\frac{\al^3}{(1-\al)^3}\right)<-M.$$
Thus $w\in C(\overline{\Om})$ solves
$$\Df w\le h(x)-\lam a(x) w^3,\;w>0,\;\mbox{in $\Om$, and $w=m_1\ge b$ on $\p\Om.$}$$

(ii) We now construct a sub-solution $v\in C(\overline{\Om})$ that satisfies
$$\Df v+\lam a(x) v ^3\ge  M,\;v<0,\;\mbox{in $\Om$, and $v=\ell_1$ on $\p\Om$.}$$
If we take $v=(\ell_1/m_1)w$, where $w$ is as in part (i), we obtain  that
$$\Df v+\lam a(x) v^3> \frac{M |\ell_1|^3}{m_1^3}\ge h(x),\;v<0,\;\mbox{in $\Om$, and $v=\ell_1\le b$ on $\p\Om$.}$$
Invoking Theorem \ref{exst0}, we obtain the existence of a solution $u\in C(\overline{\Om}),\;v\le u\le w,$ to (\ref{exst7}). 
$\Box$
\vsp
We conclude this section with a result about distance estimates regarding how close the points of a level set, of any positive solution $u$ of (\ref{Frd}), are to the boundary $\p\Om$. Define
$$F(t)=\int^1_t \; \frac{1}{(1-s^4)^{1/4}}\;ds,\;\;0\le t\le 1.$$

\begin{lem}\label{dst} Suppose that $a(x)\in C(\Om)\cap L^{\infty}(\Om),\;a(x)>0$, $\lam>0$ and $\dl\ge 0$. Let $u\in C(\overline{\Om}),\;u>0,$  solve the problem
\eqRef{dst1}
\Df u+\lam a(x) u^3=0,\;\mbox{in $\Om$, and $u=\dl$ on $\p\Om$}.
\ee
Set $\nu=\inf_{\Om}a(x)$ and $d(x)=$dist$(x,\p\Om),\;x\in \Om$.
It follows that 
$$d(x)\le \frac{F(\dl/u(x))}{(\lam \nu)^{1/4}}\le \frac{F(0)}{(\lam \nu)^{1/4}}.$$
If $m=\sup_{\Om} u$ and $z\in \Om$ is such that $u(z)=m$, then $d(z)\le F(\dl/m)/(\lam \nu)^{1/4}.$ 
\end{lem}

\NI{\bf Proof.}  First notice that the integral $F(0)<\infty.$
Let $x\in \Om$. Set $d=d(x)$ and consider the ball $B_{d}(x).$ For $0\le r\le d$, define $m(r)=\inf_{B_r(x)}u.$ 
Since $u$ is infinity superharmonic,  $m(r)=\inf_{\p B_r(x)}u$, $m(r)$ is concave and is decreasing.  Also $m(0)=u(x)$ and 
$m(d)=\dl$.

For $y\in B_d(x)$, set $r=|x-y|$. Let $w(y)=w(r)\in C(\overline{B_d(x)})$ be  defined as 
\eqRef{eqtn1}
w(r)=w(0)-(3\lam \nu)^{1/3}\int_0^r \left(\int_0^t m(s)^3\;ds\right)^{1/3}\;dt
\ee
Here $w(0)$ is so chosen that $w(d)=\dl$. Note that $w^{\prime}(0)=0$.
Using (\ref{radl}), one can show that $w$ is a viscosity solution to
$$\Df w(y)+\lam \nu m(r)^3=0,\;\mbox{in $B_{d}(x),$ and $w=\dl$ on $\p B_d(x).$ }$$
See Lemma 4.1 in \cite{BMO2} for a proof. Next, $u$ solves (\ref{dst1}), in $B_d(x)$, with $u\ge \dl$, on $\p B_d(x)$. Thus, Lemma \ref{cmp3} implies that 
$w\le u,$ and $w(r)\le m(r),$ in $B_d(x)$. Thus,
$$(w^{\prime}(r))^2w^{\prime\prime}(r)+\lam \nu w^3\le 0,\;\mbox{in $B_d(x)$, and $w(d)=\dl$.}$$
Noting that $w^{\prime}(r)\le 0$ and $w(r)>0,$ and multiplying both sides by $w^{\prime}(r)$, an integration leads to
$$(\lam \nu)^{1/4}d\le \int^{w(0)}_\dl \frac{ds}{\left(w(0)^4-s^4\right)^{1/4}}\le \int^1_{(\dl/u(x))}\frac{ds}{(1-s^4)^{1/4}}.$$
The conclusion of the lemma holds. $\Box$
\vsp
\section{\bf Existence of the first eigenvalue and the first eigenfunction}

In this section, we will show that $\lamo$, defined in (\ref{set1}), is the first eigenvalue of $\Df $ on $\Om$. The proof will also provide us with the existence of a first  
eigenfunction which turns out to be positive. As was shown in Lemma \ref{nst}, solutions to (\ref{egpb}), for $\lam<\lamo$, are the zero-solutions. Thus $\lamo$ is the smallest value of $\lam$, in (\ref{egpb}), that supports a non-trivial solution. This section also contains some monotonicity results about the first eigenvalues of the level sets of a positive first eigenfunction on $\Om$.  

In this section, we will always take $\Om\subset \IR^n$ to be a bounded domain. For a better exposition, recall (\ref{Frd}), (\ref{set}) and (\ref{set1}). In Section 3, we showed that if
$a(x)\in C(\Om)\cap L^{\infty}(\Om),\;a(x)>0,\; \dl>0 \;\mbox{and}\; 0\le \lam<\lamo,$ then there exists a positive solution $u\in C(\overline{\Om})$ to 
\eqRef{sec40}
\Df u+\lam a(x) u^3=0,\;\mbox{in $\Om$, and $u=\dl$ on $\p\Om$.}
\ee
Moreover, by Remark \ref{unq}, $u$ is unique. We recall Remark \ref{set3}, where it is shown that if $\{\lam_k\}_{k=1}^{\infty},\;\lam_k\in S,$ is an increasing sequence and if $u_k$ is the positive solution to (\ref{sec40}) corresponding to $\lam_k$, then $u_{k+1}\ge u_k$ in $\Om$. We record this fact in
\eqRef{sec42}
\mbox{$u_k,\;k=1,2,\cdots,$ is an increasing sequence.}
\ee 
We now prove the main result of this section. Also see \cite{J}.
\vsp
\begin{thm}\label{eigv} 
Let $\Om\subset \mathbb{R}^n,\;n\ge 2,$ be a bounded domain, and $a(x)\in C(\Om)\cap L^{\infty}(\Om)$ with $a(x)> 0$. Let $S$ be as defined in (\ref{set}) and $\lamo=\sup S$. Then there is a solution $v\in C(\overline{\Om}),\;v>0$ to the eigenvalue problem
$$\Df v+\lamo a(x) v^3=0\;\;\;\mbox{in $\Om$, and $v=0$ on $\p\Om$}.$$ 
\end{thm}
\NI{\bf Proof:} For $k=1,2,\cdots$, let $\lam_k \in S$, be  an increasing sequence with $\lim_{k\uparrow \infty}\lam_k=\lamo$. Fix $\dl>0$ and let $u_k>0$ solve the problem
\eqRef{eigvk}
\Df u_k+\lam_k a(x) u_k^3=0, \;\mbox{in $\Om$, and $u_k=\dl$ on $\p\Om$.}
\end{equation}
Set $m_k=\sup_{\Om} u_k$, it follows from (\ref{sec42}) that $m_k$ is increasing. We claim that
\eqRef{sec43}
\lim_{k\rightarrow \infty} m_k=\infty.
\ee
We provide a lower bound for $m_k$ by using Remark \ref{set21} and Lemma \ref{up}. By Lemma \ref{up}, for each $k=1,2,\cdots,$ there is a $\hat{u}_k>0$ such that
$$\Df \hat{u}_k+(\lam_k+\vep)a(x) \hat{u}_k^3= 0,\;\mbox{in $\Om$, and $\hat{u}_k=\dl$ on $\p\Om$,}$$  
where $0<\vep< \lam_k(\dl/m_k)^3.$ We claim that $\lamo-\lam_k\ge\lam_k(\dl/m_k)^3.$ If this were false then by taking $\vep=\lamo-\lam_k$ in Lemma \ref{up}, we would obtain a positive solution to 
$\Df \eta+\lamo a(x) \eta^3= 0,$ in $\Om$, and $\eta=\dl$ on $\p\Om.$  This would imply that $\lamo< \sup S,$ this contradicts the definition of $\lamo.$ In other words, the claim holds and
$$m_k\ge \dl \left(\frac{\lam_k}{\lamo-\lam_k}\right)^{1/3}.$$
Thus (\ref{sec43}) holds. 
 
Next, define $v_k=u_k/m_k$. Then $\sup v_k=1$ and
\eqRef{eigv0}
\Df v_k+\lam_ka(x) v_k^3=0,\;\mbox{in $\Om$, and $v_k=\dl/m_k$.} 
\end{equation}
As $v_k$'s are uniformly bounded, by Lemma \ref{lpc}, they are uniformly locally Lipschitz continuous. There is a subsequence, which we continue to denote by $\{v_k\},$ that converges locally uniformly to some function $v\in C(\Om)$ such that $v\ge 0.$ By Lemma 5.1 in \cite{BMO1}, it follows that $v$ solves 
\eqRef{eigv2}
\Df v+\lamo a(x) v^3=0,\;v\ge 0,\;\mbox{in $\Om$, and $\sup_{\Om} v=1.$}
\end{equation} 
In order to show that $v>0$ and $v\in C(\overline{\Om})$, we will employ an upper bound and a lower bound.

We first construct an upper bound. Set $\mu=\sup_{\Om}a(x)$, and let $\eta\in C(\overline{\Om})$ solve the problem
$$\Df \eta=-2\lamo\mu,\;\mbox{in $\Om$, and $\eta=0$.}$$ 
The existence of $\eta$ follows from \cite{BMO1,BMO2,LW1}. Also,  the function $\eta+\dl/m_k$ solves the same differential equation with $\dl/m_k$ as the boundary data. 
Since (\ref{eigv0}) implies that $\dl/m_k\le v_k\le 1$, it is easy to see that $2\lamo\mu\ge \lamo a(x) v_k^3$. It follows from (\ref{eigv0}) and Lemma \ref{cmp2}, 
that  $\dl/m_k\le v_k\le \eta+\dl/m_k,\;k=1,2,\cdots.$ Thus $0\le v\le \eta$, in particular, $v=0$ on $\p\Om$ and 
$v\in C(\overline{\Om})$. In order to show that $v>0$ in $\Om$, we construct a lower bound. Since $v_k$'s are continuous in $\Om$ and $\sup_{\Om}v_k=1$, there is a point $x_k\in \Om$ such that
$v_k(x_k)=1$. We may now find a subsequence of $v_k$ and $x_k$( which we continue to call them as $v_k$ and $x_k$)  with $x_k\rightarrow x$. Since $v$ is small near $\p\Om$, it follows that $x\in \Om$ and $v(x)=1$. Let $h\in C(\overline{\Om}\setminus\{x\})$ solve
$$\Df h=0\;\;\mbox{in $\Om\setminus\{x\}$, with $h(x)=2/3$ and $h=0$ on $\p\Om$.}$$
By Lemma \ref{cmp3} and (\ref{eigv0}), $0<h\le v_k$, for large $k$. Thus, $0<h\le v \le \eta$. The conclusion of the theorem follows.\quad $\Box$  
\vsp
From hereon we will refer to $\lamo$ as the first eigenvalue of the infinity-Laplacian and a non-trivial solution $u\in C(\overline{\Om})$ to the problem
\eqRef{eignfn}
\Df u+\lamo a(x) u^3=0,\;\mbox{in $\Om$, and $u=0$ on $\p\Om$,}
\ee
as a first eigenfunction. As is clear from Theorem \ref{eigv}, an eigenfunction, having one sign in $\Om$, exists. In the rest of this section, we will derive some properties of $\lamo$.
We  start with an observation about domain monotonicity of the first eigenvalue. 
\vsp
\begin{rem}\label{mon}
\NI In (\ref{set}), let us write $S=S(\Om)$. Suppose that 
$\Om^{\prime}\subset \Om$ is a sub-domain. If $\lam>0$ is such that there is a function $v\in C(\overline{\Om})$ that solves 
$$\Df v+\lam a(x) v^3=0,\;\;v>0,\;\mbox{in $\Om$, and $v=\dl$ on $\p\Om$,}$$
then $v$ also solves the same equation in $\Om^{\prime}$ with $v\ge \dl$ on $\p\Om'$. Thus
$S(\Om)\subset S(\Om^{\prime})$ and $\lamo\le \lam_{\Om^{\prime}}$.

Suppose that $\Om'$ is compactly contained in $\Om$ and $u>0$ solves (\ref{eignfn}), see Theorem \ref{eigv}. If we set $\tht=\inf_{\Om^{\prime}} u,$ then $\tht>0$.
Since one can use $u$ as a super-solution and the function $v=\tht$ as a sub-solution of (\ref{eignfn}), Theorem \ref{exst0} provides us with a positive solution $w\in C(\overline{\Om^{\prime} } )$ to 
the problem
$$\Df w+\lamo a(x) w^3=0,\;\mbox{in $\Om^{\prime}$, and $w=\tht$ on $\p\Om^{\prime}$.}$$
By Lemma \ref{up}, we can find an $\vep>0$ and a function $w\in C(\overline{\Om^{\prime}})$ that solves
$$\Df \bar{w}+(\lamo+\vep) a(x) \bar{w}^3=0,\;\bar{w}>0,\;\;\mbox{in $\Om^{\prime}$, and $\bar{w}=\tht$ on $\p\Om^{\prime}$.}$$
By the definition of the set $S$, we see that $\lam_{\Om^{\prime}}\ge \lamo+\vep>\lamo.$ We have thus strict domain monotonicity in case $\Om^{\prime}$ is compactly contained in $\Om$. However, in general, there is no strict domain monotonicity, see Lemma \ref{annb} in Section 5.\quad $\Box$  
\end{rem}
\vsp
\begin{rem}\label{eigv1}
We also observe that if $u\in C(\overline{\Om}),\;u\not=0,$ solves 
$$\Df u+\lam a(x) u^3=0,\;\;\mbox{in $\Om$},\;\mbox{and $u=0$ on $\p\Om$,}$$
then $\lam\ge \lam_{\Om}.$ This can be seen as follows. Firstly, by Lemma \ref{ap}, $\lam>0$. Next, if $\lam<\lamo$, then Remark \ref{set21} and Lemma \ref{nst} would imply that $u=0$ in $\Om$. Thus the claim holds. $\Box$ 
\end{rem}

In the next lemma, we make an observation related to Remark \ref{mon}. This addresses the monotonicity property of the first eigenvalue of a level set of an eigenfunction.
\begin{lem}\label{mnct}
Let $\Om\subset \mathbb{R}^n,\;n\ge 2$, be a bounded domain, and $a(x)\in C(\Om)\cap L^{\infty}(\Om)$ with $a(x)> 0$. Let $u\in C(\overline{\Om}),\;u >0,$ and $\sup_{\Om} u=1$ be a first eigenfunction, that is,  
$$\Df u+\lamo a(x) u^3=0,\;\mbox{in $\Om$, and $u=0$ on $\p\Om$}.$$ 
For $0<t\le 1$, set $\Om_t=\{x\in \Om:\;u(x)>t\}$. Then $\lam_{\Om_t}$ is increasing and $\lim_{t\uparrow 1} \lam_{\Om_t}=\infty$.  
\end{lem} 
\NI{\bf Proof:}  First note that by Remark \ref{mon}, $\lamo<\lam_{\Om_t}<\lam_{\Om_s}$, for $0<t<s<1$, and
\eqRef{mnct1}
\Df u+\lamo a(x) u^3=0\;\;\;\mbox{in $\Om_t$, and $u=t$ on $\p\Om_t,\;\forall\;0\le t<1.$}.
\ee
For notational ease, call $\lam_t=\lam_{\Om_t}$. Now, for any fixed $0<\al<1$ and $0<t<1$, and, for any $t\le s\le 1$,  
\eqRef{mnct01}
s^3-(s-\al t)^3=3\al st(s-\al t)+\al^3 t^3\ge \al^3t^3.
\ee
Take $0<\al<1$  and $\vep>0$ to be chosen later. As done in Lemma \ref{up}, we write $w=u-\al t$ and obtain
\eqRef{mnct2}
\Df w+(\lamo+\vep)a(x) w^3=a(x)\left((\lamo+\vep) (u-\al t)^3-\lamo u^3\right),\;\;\mbox{in $\Om_t$,}
\ee
with $w=(1-\al)t$ on $\p\Om_t$. Rearranging the right side we obtain that 
$$a(x)\left\{(\lamo+\vep) (u-\al t)^3-\lamo u^3\right\}= a(x)\left[\vep (u-\al t)^3-\lamo\left\{u^3-(u-\al t)^3\right\}\right].  $$
Using (\ref{mnct01}) and $t<u\le 1$ in (\ref{mnct2}), we conclude  
$$\Df w+(\lamo+\vep)a(x) w^3\le a(x)\left(\vep(1-\al t)^3 -\lamo \al^3t^3\right),$$
For $0<\tht<1$, select
$$\vep_{\tht}=\tht \frac{\al^3 t^3 \lamo}{(1-\al t)^3},$$
to obtain 
$$\Df w+(\lamo+\vep_{\tht})a(x) w^3\le 0,\;\mbox{in $\Om$, and $w=t(1-\al)>0$ on $\p \Om_t.$}$$
By Remark \ref{set21},
$$\lam_t \ge \lamo\left(1+\tht  \frac{\al^3 t^3}{(1-\al t)^3}\right).$$
By Remark \ref{mon},
$$\lim_{t\uparrow 1} \lam_t\ge \lam_t\ge \lamo\left(1+\tht  \frac{\al^3 t^3}{(1-\al t)^3}\right).$$
The inequality holds for any $0<\al<1$ and $0<t<1$, hence the claim. \quad $\Box$ 
\vsp
We make a related observation regarding $\lamo$. In the previous lemma, we discussed the limit $\lim_{t\uparrow 1} \lam_t$. In the next lemma we study the limit
$\lim_{t\downarrow 0} \lam_t.$

\begin{lem}\label{lvl}
Suppose that $a(x)\in C(\Om)\cap L^{\infty}(\Om)$ with $a(x)>0.$ Let $T$ be the set of all $\lam$'s such that $\lam\ge \lamo$ and the following problem
$$\Df v+\lam a(x) v^3=0,\;\mbox{in $\Om$, and $v=0$ on $\p\Om$}.$$
has a positive solution $v\in C(\overline{\Om}).$ Let $u>0$ be an eigenfunction corresponding to $\lamo$. Assume that $\sup_{\Om} u=1$.
For $0<t<1$, define $\Om_{t}=\{x:\;u(x)>t\}$ and  $\lam_t=\lam_{\Om_t}.$ Then
$$\lamo=\inf T\le \sup T=\inf_t\lam_t  =\lim_{t\downarrow 0}\lam_t.   $$
In particular, $T$ is a singleton set if and only if $\lamo=\lim_{t\downarrow 0}\lam_t.$
\end{lem}

\NI{\bf Proof:} Firstly, $\sup T<\infty$, by Theorem \ref{eqn1}. If $0<t<s<1$ then $\Om_s\subset \Om_t$, and  by Remark \ref{mon}, $\lam_t\le \lam_s$ and 
$\lim_{t\downarrow 0}\lam_t=\inf_t \lam_t$.
Our goal is to show that $\lam_t \ge \sup T,$  for all $0<t<1.$ Suppose not. Let $\lam\in T$ be such that $\lam>\lam_t$, for some $0<t<1$. By the definition of $T$, there is a function $v$ that solves
$$\Df v+\lam a(x) v^3=0,\;v>0,\;\;\mbox{in $\Om$, with $v=0$ on $\p\Om$}.$$
Since $\Om_t$ is compactly contained in $\Om$, $\inf_{\p\Om_t}v>0$. Next, let $w>0$ be a first eigenfunction on $\Om_t,$ that is,
$$\Df w+\lam_t a(x) w^3=0,\;\mbox{in $\Om_t$, and $w=0$ on $\p\Om_t$.}$$ 
Since $\lam>\lam_t$, applying Lemma \ref{cmp1} to $v$ and $w$ in $\Om_t$, we obtain $0\le (w/v)\le \sup_{\p\Om_t}(w/v)=0$, a contradiction. Hence, $\lam_t\ge \sup T$, for all $0<t<1$. By Theorem \ref{eigv},
$$\lamo=\inf T\le \sup T\le \inf_t \lam_{t}.$$

Next, we show that $\sup T= \inf_{t}\lam_{t}.$ To see this,   for $0<t<t_0$, $t_0$ small, consider the family of first eigenfunctions  $w_t$ that solve
$$\Df w_{t}+\lam_{t} a(x) w_{t}^3=0,\;w_{t}> 0,\;\mbox{in $\Om_t$, and $w_t=0$ on $\p\Om_t$}.$$
Scale $w_t$ such that $\sup w_t=1.$ Calling $\lam_0=\inf_t \lam_t$ and arguing as in Theorem \ref{eigv} (see Lemma 5.1 in \cite{BMO1}), we obtain a convergent subsequence 
 $\{w_{t_l}\}_{l=1}^{\infty}$ (with $t_l\downarrow 0$) of $\{w_t\}_{t<t_0}$ and a function $w_0\in C(\Om)$ such that $\lim_{t_l\rightarrow \infty} w_{t_l}=w_0$ with $\sup w_0=1$. Also,
$$\Df w_0 +\lam_0 a(x) w_0^3=0,\;w_0\ge 0,\;\mbox{in $\Om$, and $w_0\ge 0$ on $\p\Om$}.$$
To show that $w_0\in C(\overline{\Om})$ and $w_0=0$ on $\p\Om$, we employ an upper bound similar to that in Theorem \ref{eigv}. Set $\mu=\sup_{\Om} a$ and let $\eta\in C(\overline{\Om})$ be the solution to
$$\Df \eta=-2\lam_{t_0}\mu,\;\;\mbox{in $\Om$ and $\eta=0$ on $\p \Om$.}$$
Since for any $0<t<t_0$, $\Om_t\subset \Om,\;\lam_t\le \lam_{t_0},\; 0\le w_t\le 1$ and $\eta>0$ in $\Om_t$, Lemma \ref{cmp3} implies that $w_t\le \eta$ in $\Om_t$. Thus $0\le w_0\le\eta$ in $\Om$, and thus, $w_0\in C(\overline{\Om})$ and $w_0=0$ on $\p\Om$. 

We now prove that $w_0>0$ in $\Om$. Let $x_l\in \Om_{t_l},\;l=1,2,\cdots,$ be such that $w_{t_l}(x_l)=1$. Then for some $p\in \overline{\Om}$, $x_l\rightarrow p$ as $l\rightarrow \infty$ (choose a subsequence, if needed). Since $w_{t_l}\le \eta$, it follows that $p\in \Om.$ Hence, $w_{t_l}(p) >1/2$, for $t_l$ close to $0$. Take $s$, close to $0$, such that $p\in \Om_s$(any $s<u(p)$ will do). 
We take $\zeta$ to be a positive infinity harmonic function in $\Om_s\setminus \{p\}$ with $\zeta(p)=1/2$ and $\zeta=0$ on $\p\Om_s$. Since $w_{t_l}$ is positive and infinity super-harmonic in $\Om_{t_l}$ and $\Om_s\subset \Om_{t_l}$, for $0<t_l<s$, Lemma \ref{cmp2} implies that 
$w_{t_l}\ge \zeta$ in $\Om_s$. Thus $w_0\ge \zeta>0$ in $\Om_s$. In particular, $w_0>0$ in $\Om_s$, for any $s$ close to $0$. Since $\Om_s$ exhausts $\Om$ as $s$ decreases to $0$, we have that $w_0>0$ in $\Om$.

Thus  $\inf_t\lam_t=\sup T.$ The claim holds. \quad $\Box$ 
\vsp
\begin{rem}\label{cm} Let  the function $u\in C(\overline{\Om})$, the sets $\Om_t$, the eigenvalues $\lam_t,\;0<t<1,$ and $T$ be as in the statement of Lemma \ref{lvl}. We claim that the set $T$ is either a singleton set or the interval $[\lamo,\;\sup T].$ Set $\lam^T=\sup T$, and assume that $T$ is not a singleton set. Choose $\vep>0$ such that $\lam^T-\vep>\lamo$. 
 Fix $\dl>0$, and for each $0<t<1$,
consider the family of problems
$$\Df v_{t}+(\lam^T-\vep) a(x) v_{t}^3=0,\;v_{t}>0,\;\mbox{in $\Om_{t}$, with $v_{t}=\dl$ on $\p\Om_{t}$}.$$
By Lemma \ref{lvl}, $\lam_t>\lam^T-\vep$. Hence, Theorem \ref{exst3} (also see Remark \ref{set21}) implies that the above has a unique solution $v_t\in C(\overline{\Om_t})$, $v_t>\dl$, for every $0<t<1.$ If $0<t_1<t_2<1,$ then $\Om_{t_2}\subset \Om_{t_1}$ and $\lam_{t_1}<\lam_{t_2}$, and we conclude from Lemma \ref{ord3} that $v_{t_2}\le v_{t_1}$, in $\Om_{t_2}.$
Call $m_{t}=\sup_{\Om_t} v_{t}$, then $m_{t}$ increases as $t$ decreases. We claim that $\lim_{t\downarrow 0}m_t=\infty$. To see this, first we employ Lemma \ref{cmp1},  noting that $\lam_t>\lamo$,
to observe that $\sup_{\Om_t}(u/v_t)=\sup_{\p\Om_t}(u/v_t)=t/\dl.$  If $\sup_t m_t<\infty$ then it follows that $u\le (t m_t)/\dl,$ in $\Om_t$.
Letting $t$ decrease to $0$, we get $u=0$ in $\Om$. This is a contradiction and the claim holds.

Define $w_t=v_t/m_t$, in $\Om_t$. Noting that $\sup_{\Om_t}w_t=1$ and arguing as in Theorem \ref{eigv} and Lemma \ref{lvl} (see Lemma 5.1 in \cite{BMO1}), one can find a convergent subsequence $\{w_{t_l}\}$ of $\{w_t\}$(with $t_l\rightarrow 0$) and $w\in C(\overline{\Om})$ such that $\lim_{t_l\rightarrow 0}w_{t_l}\rightarrow w$. Moreover, 
$$\Df w+(\lam^T-\vep) a(x) w^3=0,\;w>0,\;\;\mbox{in $\Om$, with $w=0$ on $\p\Om$}.$$
This proves our assertion. $\Box$
\end{rem}
\vsp
\section{\bf Additional results on some special domains}

In Sections 5 and 6, we will discuss some results regarding the first eigenvalue problem on some special domains. The present section contains a discussion related to
the eigenvalue problem (\ref{eignfn}) on $C^2$ domains and on star-shaped domains. If $\lam_t$ and $T$ are as in the statement of Lemma \ref{lvl}, we will show that $T$ is a singleton set when $\Om$ is a $C^2$ domain, in other words, $\lim_{t\downarrow 0}\lam_t=\lamo$, see Remark \ref{sng}.

We begin this section by proving that the eigenfunctions corresponding to higher eigenvalues change sign. This fact is well-known in the context of elliptic operators on general domains. 
We provide a proof in this context for $C^2$ domains and star-shaped domains. In this context, recall the result in Theorem \ref{eqn1} that holds on any bounded domain.

\begin{lem}\label{sgn} Let $\Om\subset \IR^n$ be a bounded domain. Suppose that either $\Om$ has $C^2$ boundaries or is star-shaped. We assume that 
(i) $a(x)\in C(\overline{\Om})$ 
and, $\inf_{\Om} a(x)>0,$ if $\Om$ is star-shaped, and (ii) that $a(x)\in C(\Om)\cap L^{\infty}(\Om),\;a(x)>0,$ if $\Om$ has a $C^2$ boundary.
Let $\lam>\lamo$ and $v\in C(\overline{\Om})$ be such that 
 \eqRef{sgn1}
 \Df v+\lam a(x) v^3=0,\;\mbox{in $\Om$, $\sup_{\Om} v=1,$ and $v=0$ on $\p\Om$.}
 \ee
Then $v$ changes sign in $\Om$.  
\end{lem}
\NI{\bf Proof.} We start with the case when $\Om$ is a star-shaped domain. Without any loss of generality, we may assume that $\Om$ is star-shaped 
with respect  to the origin $o$. Suppose  that $v>0$ in $\Om$. We scale $v$ as follows.
For $0<t<\infty$, set $y=tx$, $w_t(y)=v(x)$ and $\Om_t=\{tx:\;x\in \Om\}$. Note that 
$\Om_s\subset \Om\subset \Om_t,\;0<s<1<t.$
A simple calculation leads to  
$$\Df w_t+\frac{\lam}{t^4}a(y/t)w_t^3=0,\;\mbox{in $\Om_t$, and $w_t=0$ on $\p\Om_t$.}$$
Taking $t>1$, close to $1$, and using the uniform continuity of $a$, we have
$$\lamo a(y) \le \frac{\lam}{ t^4}a(y/t),\;y\in \Om.$$ 
Hence,
$$\Df w_t+\lamo a(y) w_t^3\le 0,\;\mbox{in $\Om$, and $\inf_{\p\Om}w_t>0$.}$$
This contradicts the definition of $\lamo$, see Remark \ref{set21} and Lemma \ref{up}. The claim holds.

We now prove the lemma when $\Om$ is $C^2$. We achieve this in six steps. We assume that $v>0$ in (\ref{sgn1}).

\NI{\bf Step 1:} By Theorem \ref{eigv}, one can find an eigenfunction $u>0$ such that
\eqRef{sgn2}
\Df u+\lamo a(x) u^3=0,\;\mbox{in $\Om$, $\;\sup_{\Om} u=2,$ and $u=0$ on $\p\Om$.}
\ee

\NI{\bf Step 2:} We construct two auxiliary functions. Set $\s=3^{4/3}/4$, and consider the ball $B_{R}(o)$, for $R>0.$ Take $m>0$, define
\eqRef{sgn3}
\psi(x)=\psi(m,R,|x|)=c|x|-b|x|^{4/3},\;\;\;x\in B_R(o),
\ee
where $c=(1/R)+(8\s/3)R^{1/3}m^{1/3},$ and $b=\s m^{1/3}.$
Using (\ref{radl}), for $x\ne o$, we have
\ben
\Df \psi&=&\left(-\frac{4b}{9}\right)\left(c-\frac{4b|x|^{1/3}}{3}\right)^2 |x|^{-2/3}\\
&=&\left(-\frac{4\s}{9}\right)\left[m^{1/3}|x|^{-2/3} \left(c^2-\frac{8 cb|x|^{1/3}}{3}+\frac{16 b^2|x|^{2/3}}{9}\right)\right]\\
&=&\left(-\frac{4\s}{9}\right)\left[m^{1/3}|x|^{-2/3}\left(c^2-\frac{8cb|x|^{1/3}}{3}\right)+\frac{16m^{1/3}b^2}{9}\right]\\
&=&\left[\left(-\frac{4\s}{9}\right)c m^{1/3}|x|^{-2/3}\left(c-\frac{8\s m^{1/3}|x|^{1/3}}{3}\right)\right]-m\\
&=&\left(-\frac{4\s}{9}\right)c m^{1/3}|x|^{-2/3}\left(\frac{1}{R}+\frac{8\s m^{1/3}}{3}\left(R^{1/3}-|x|^{1/3}\right)\right)-m\le -m.
\een
We record this and other useful facts for $\psi$, see (\ref{sgn3}), 
\eqRef{sgn4}
(i)\;\psi(o)=0,\;\;(ii)\; \psi(R)>1,\;\mbox{and}\;\;(iii)\;\Df \psi(x)\le -m,\;\;\psi(x)>0,\;\;x\in B_{R}(o)\setminus\{o\}.
\ee
 For $\ell>0$, define
$$\eta(x)=\eta(\ell,R,|x|)=\ell \left(1-\frac{|x|}{R}\right),\;\;\forall\;x\in B_{R}(o).$$  
We note also the following for future reference.
\eqRef{sgn5}
(i)\;\;\eta(R)=0,\;\;(ii)\;\;\eta(o)=\ell,\;\;\mbox{and}\;\;(iii)\;\;\Df \eta(x)=0,\;\;x\in B_{R}(o)\setminus\{o\}.
\ee

 We introduce additional notations that will be used in Steps $3,\;4$ and $5$. Being a $C^2$ domain, $\Om$ satisfies an uniform interior ball condition at every point of $\p\Om.$ Let $2\rho$ denote the radius of the optimal ball.  For every $z\in \p\Om$, let $\nu(z)$ denote the unit inward pointing normal. Then the ball $B_{2\rho}(z+2\rho\nu(z))\subset \Om$ and $z\in \p\Om\cap \p B_{2\rho}(z+2\rho\nu(z)).$ For every $z\in \p\Om$, set $y=z+\rho \nu(z).$  

\NI {\bf Step 3:} For every $z\in \p\Om$, define
\eqRef{sgn6}
\Om_*=\Om\setminus \left(\cup_{z\in \p\Om}B_{\rho/2}(z+\rho\nu(z)/2)\right).
\ee
Also, set
\eqRef{sgn7}
\ell_u=\inf_{\Om_*} u\;\;\;\;\mbox{and}\;\;\;\ell_v=\inf_{\Om_*} v,
\ee
where $u$ is as in Step 1 and $v$ is as in (\ref{sgn1}). \\
\NI {\bf Step 4:} We work in the balls $B_{\rho}(y)$ and $B_{2\rho}(z)$. Here, $B_{\rho}(y)\subset \Om\cap B_{2\rho}(z).$ 
We recall the constructions in Step 2, (\ref{sgn3})-(\ref{sgn5}) and (\ref{sgn7}). Let $\mu=\sup_{\Om} a(x).$ Recalling Step 1, take $m_u=\lamo \mu$ and $m_v=\lam \mu$. For each fixed
$z\in \p\Om$, set in (\ref{sgn3}),
\eqRef{sgn71}
\psi_u(x)=\psi(m_u, 2\rho, |x-z|)\;\mbox{and}\;\;\psi_v(x)=\psi(m_v,2\rho,|x-z|),\;\;x\in B_{2\rh}(z).
\ee
Next, in Step 2, take
\eqRef{sgn8}
\eta_u(x)=\eta(\ell_u, \rho,|x-y|)\;\;\mbox{and}\;\;\eta_v(x)=\eta(\ell_v,\rho,|x-y|),\;\;x\in B_{\rh}(y).
\ee
We also note that if $x\in B_{\rh}(y)$ and lies on the segment $yz$, then 
\eqRef{sgn9}
\eta_u (x)=\frac{\ell_u |x-z|}{\rho}\;\;\mbox{and}\;\; \eta_v(x)=\frac{\ell_v|x-z|}{\rho}.
\ee
\NI{\bf Step 5:} We claim that for each $z\in \p\Om$ and $x\in B_{\rho}(y)$
\eqRef{sgn10}
\eta_u(x)\le u(x)\le 2\psi_u(x),\;\;\;\mbox{and}\;\;\;\eta_v(x)\le v(x)\le \psi_v(x).
\ee
We present details for $u$, the proof for $v$ will follow analogously. We apply the properties of $\psi_u$ from (\ref{sgn4}) in $B_{2\rho}(z) \cap \Om$, call $w=2 \psi_u$. 
Using Step 1, (\ref{sgn4}) and (\ref{sgn71}), we see that 
$$\Df w\le-8\lamo\mu,\;\;\mbox{and}\;\;\Df u\ge -8 \lamo\mu,\;\;\mbox{in $B_{2\rho}(z) \cap \Om$.}$$
From (\ref{sgn2}), (\ref{sgn4}) (ii) and (iii), we see that $w(x)>2\ge u(x),\; x\in \p B_{2\rho}(z)\cap \Om$, and $w\ge u$ on $\p\Om\cap B_{2\rho}(z)$.
The comparison principles in the Lemmas \ref{cmp2} and \ref{cmp3} yield that $u\le 2\psi_u$, in $B_{2\rho}(z)\cap \Om$. To show that $\eta_u\le u$, in $B_{\rho}(y)$, we note
$$\Df \eta_u=0,\;\;\mbox{and}\;\;\Df u\le 0,\;\;\mbox{in $B_{\rho}(y)\setminus\{y\}$.}$$
Using (\ref{sgn5})-(\ref{sgn7}) and (\ref{sgn9}), we have that $\eta_u(x)\le u(x),\;x\in \p B_{\rho}(y)$ and $\ell_u=\eta_u(y)\le u(y)$. 
Thus, Lemma \ref{cmp2} implies that $\eta_u\le u$ in
$\overline{B_{\rho}(y)}$. Thus (\ref{sgn10}) holds. 

If $x\in \Om\setminus \overline{\Om_*}$ (see (\ref{sgn6})), then one can find a closest point $z\in \p\Om$, such that $x\in B_{\rho}(y)$, where $y=z+\rho \nu(z)$. As a result, we have
$$\frac{\eta_u(x)}{\psi_v(x)}\le \frac{u(x)}{v(x)}\le \frac{2\psi_u(x)}{\eta_v(x)}.$$
Next, we observe that $x$ lies on the segment $yz$. From Step 2 and (\ref{sgn9}), we conclude that there are positive constants $k_1$, $k_2$ and $d$, depending only on $\ell_u,\;\ell_v,\;\lam,\;\lamo,\;\mu$ and $\rho$, such that
\eqRef{sgn11}
k_1\le \frac{u(x)}{v(x)}\le k_2,\;\;\mbox{for every $x\in \Om$ with dist$(x,\;\p\Om)<d.$}
\ee
\NI{\bf Step 6:}  We recall (\ref{sgn1}), (\ref{sgn2}), (\ref{sgn11}) and Lemma \ref{cmp}. Choose $1<\tau<(\lam/\lamo)^{1/3}.$ Since $u-v=0$ on $\p\Om$, $\sup_{\Om} u=2$ and $\sup_{\Om} v=1$, the function $u-v$ will assume a positive maximum in $\Om$. We will show that this leads to a contradiction thus proving the lemma. 

Since $\sup_{\Om}(u-v)>\sup_{\p\Om}(u-v)$, by Lemma \ref{cmp}, there is a point $x_1\in\Om$, where $u-v$ takes its supremum and
$(\lam/\lamo)^{1/3}v(x_1)\le u(x_1).$ As $(u-\tau v)(x_1)>0$ and $(u-\tau v)=0$ on $\p\Om$, the function $u-\tau v$ has a positive maximum in $\Om$. An application of Lemma \ref{cmp} to $u$ and $\tau v$ yields that there is an $x_2\in \Om$ such that
$$\sup_{\Om}(u-\tau v)=(u-\tau v)(x_2)>0,\;\;\;\mbox{and}\;\;\;\tau\left(\frac{\lam}{\lamo}\right)^{1/3}v(x_2)\le u(x_2).$$
We iterate this argument. Suppose that we have shown for some $m=1,2,\cdots,$ that there is an $x_m\in \Om$ such that 
$$\sup_{\Om}(u-\tau^{m-1}v)=(u-\tau^{m-1} v)(x_m)>0,\;\;\mbox{and}\;\;\;\tau^{m-1}\left(\frac{\lam}{\lamo}\right)^{1/3}v(x_m)\le u(x_m).$$
Since $u-\tau^{m}v=0$ on $\p\Om$, the function $u-\tau^m v$ has a positive maximum in $\Om$. Applying Lemma \ref{cmp} to $u$ and $\tau^m v$, we see that there is an $x_{m+1}\in \Om$ such that
$$\sup_{\Om}(u-\tau^{m}v)=(u-\tau^{m} v)(x_{m+1})>0,\;\;\mbox{and}\;\;\;\tau^{m}\left(\frac{\lam}{\lamo}\right)^{1/3}v(x_{m+1})\le u(x_{m+1}).$$
Thus, we have shown that for each $m=1,2,\cdots$, there is an $x_m\in \Om$ such that 
$$u(x_m)\ge \tau^m v(x_m).$$ 
Recall that the functions $u$ and $v$ are in $C(\overline{\Om})$, $u>0,\;v>0$, in $\Om$, and $u=v=0$ on $\p\Om$. It follows that
$v(x_m)\rightarrow 0$ as $m\rightarrow \infty$. that is, $x_m$ is close to $\p\Om$ for large $m$. Combining this with (\ref{sgn11}), we obtain
$\tau^m\le k_2,$
for all values of $m$ that are large enough. This is a contradiction and the lemma holds. Incidentally, (\ref{sgn1}), (\ref{sgn2}), (\ref{sgn11}) and Lemma \ref{cmp1} lead to
$u\le k_2 v$ in $\Om$. This could have been used instead to achieve the last part of the proof.  $\Box$

\begin{rem}\label{sng} Lemma \ref{sgn} leads to the following conclusions.\\
\NI (i) Suppose that $\lam=\lamo$, and $u$ and $v$ are two positive eigenfunctions. Adapting the arguments in Step 2-5 of Lemma \ref{sgn} and applying Remark \ref{unq}, 
we have that $k_1\le u/v\le k_2,\;\mbox{in $\Om$.}$\\
\NI (ii) By Lemma \ref{lvl}, Remark \ref{cm} and Lemma \ref{sgn}, it follows that $\lim_{t\downarrow 0}\lam_{t}=\lam^T=\lamo.$ Thus $T$ is a singleton set. $\Box$
\end{rem}

\section{\bf Case of the ball}

 We now turn our attention to the case of the ball. We will take the weight function $a(x)$ to be radial. We will study the radial version of the eigenvalue problem and present some properties of the radial eigenfunction. Under the hypothesis that $a(x)$ is a constant function, we provide a description of the eigenvalues that support radial eigenfunctions and show that there are infinitely many such eigenvalues. We end the section by presenting a proof of the fact that if the weight function is a constant then the first eigenfunction has one sign and all radial first eigenfunctions are unique up to scalar multiplication. 

We begin by recalling that the existence of the first eigenvalue and a positive first eigenfunction is guaranteed by Theorem \ref{eigv}. We apply now the results of Section 3 and 4 to show that there is a first eigenfunction $u$ that is positive and radial.

For $R>0,$ let $\Om=B_R(o),$ and we take $a(x)=a(|x|)>0$. For ease of notation, we set $\lam_B=\lam_{B_R(o)}$ and $r=|x|.$ If $v(x)=v(r)$ then the radial expression for the infinity-Laplacian in (\ref{radl}) gives us  
\eqRef{rad0}
\Df v+\lam a(x) v^3=\left(\frac{dv}{dr}\right)^2\frac{d^2v}{dr^2}+\lam a(r) v(r)^3,\;\;x\in B_R(o).
\ee
 Let us also recall from Section 3 the following definition of $F(t)$ for $0\le t\le 1$, that is,
\eqRef{rad000}
F(t)=\int_t^1 \frac{ds}{(1-s^4)^{1/4}}.
\end{equation}
The ideas of the proof of Theorem \ref{rad} and Lemma \ref{annb}, that follow, are similar to those in Lemma 6.1 in \cite{BMO2}. 
\begin{thm}\label{rad}
Let $a(x)\in C(B_R(o))\cap L^{\infty}(B_R(o)),\;a(x)>0,$ and $\lam>0$. Assume that $a(x)=a(|x|)$. Let $\dl\ge 0$, and $u$ solve
\eqRef{rad1}
u(x)=u(r)=m-(3\lam)^{1/3}\int_0^r \left[ \int_0^t a(s) u(s)^3\;ds\right]^{1/3}\;dt,
\ee
where $u(o)=m>0$ is so chosen that $u(R)=\dl$. Then $u\in C(\overline{B_R(o)})$ and the following hold.\\
\NI (i) If $\lam<\lam_B$ and $\dl>0$ in (\ref{rad1}), then $u>0$, in $B$, and $u$ is the unique solution to 
\eqRef{rad2}
\Df u+\lam a(x) u^3=0,\;u>0,\;\mbox{in $B_R(o),$ and $u(R)=\dl$.}
\ee
(ii) If $\lam=\lam_B$, in (\ref{rad1}), then there is a positive function $v$ that solves (\ref{rad1}) in $B$, with $m=1$, $v(R)=0$. Moreover, $v$ is a radial first eigenfunction. \\
\NI (iii) Let $a(x)=k$ be a positive constant and $F$ be as (\ref{rad000}). Then the positive function $u$ defined by
\eqRef{rad3}
F(u(r)/m))=(\lam k)^{1/4} r,
\ee
is a radial solution to (\ref{rad2}) with $\dl\ge 0.$ We also have, $(\lam k)^{1/4}R=F(\dl/m).$  
\end{thm}    

\NI{\bf Proof:} We have broken up the proof into five steps. We take $\dl>0$. Set $\mu=\sup_{B}a(x)$ and $\nu(r)=\inf_{B_r(o)}a(x).$ 

\NI {\bf Step 1.} For any $m>\dl$, define $u$ to be the local solution to (\ref{rad1}).
By Picard's iteration, $u$ exists near $o$ and is decreasing in $r$. Since $u\in C^2$, near $o$ (except perhaps at $o$), we obtain by a differentiation that $u$ solves 
(\ref{rad2})(see (\ref{rad0})) in $r>0$, for small $r$. 
 We record a simple estimate. For small $r>0$, since, $u(r)\le u(s)\le m$, for $0\le s\le r,$ we have that
\eqRef{rad001}
\frac{(3^{4}\nu(r)\lam)^{1/3}u(r)r^{4/3}}{4}\le (3\lam)^{1/3}\int_0^r \left[ \int_0^t a(s) u(s)^3\;ds\right]^{1/3}\;dt\le \frac{(3^{4}\mu\lam)^{1/3}mr^{4/3}}{4}.
\end{equation}
\vsp
\NI{\bf Step 2.} We show that $u$ is a viscosity solution to the differential equation in (\ref{rad2}), in a neighborhood of $o$. Assume that for some $\psi\in C^2(B_R(o))$, $u-\psi$ has a local maximum at $o$, that is, $u(x)-u(o)\le \psi(x)-\psi(o)$, for $x$ near $o$. Employing (\ref{rad1}), (\ref{rad001}) and noting that $r=|x|$, we have 
$$ - \frac{(3^{4}\mu\lam)^{1/3}m|x|^{4/3}}{4}\le u(x)-u(o)\le \langle D\psi(o), x\rangle+o(|x|),\;\;\;\mbox{as}\;|x|\rightarrow 0.$$ 
Take $x=-\tht D\psi(o),\;\tht>0.$ Next, dividing both sides by $\tht$ and letting $\tht\rightarrow 0$, 
we get $D\psi(o)=0$. Hence, $\Df \psi(o)+\lam a(o) u(o)^3\ge 0$, and $u$ is a sub-solution to (\ref{rad1}). 

Suppose that $u-\psi$ has a minimum at $o$, that is
$\psi(x)-\psi(o)\le u(x)-u(o)\le 0$. Using (\ref{rad1}) and (\ref{rad001}) and arguing as above, we see that $D\psi(o)=0$. Clearly, now (\ref{rad1}) and (\ref{rad001}) lead to  
$$\frac{\langle D^2\psi(o) x, x\rangle}{2}+ o(|x|^2)\le u(x)-u(o)\le -\frac{(3^{4}\nu(r)\lam)^{1/3}u(r)|x|^{4/3}}{4},\;\;\mbox{as}\;|x|\rightarrow 0.$$  
Taking, for instance, $x=re_1$, dividing both sides by $r^2$ and then letting $r\rightarrow 0$, we see that  $D^2\psi(o)$ does not exist. Thus, $u-\psi$ can not have a minimum at $o$. Clearly, $u$ is a super-solution and, hence, a local solution to (\ref{rad2}).  
\vsp
\NI{\bf Step 3.} Steps 1 and 2 show that for any $m>\dl$, the formula in (\ref{rad1}) provides a local radial solution to (\ref{rad2}). 
By Step 1, $u$ exists near $o$ and $u$ is decreasing. Let $\vep>0$ be small. For $r>\vep$, an integration of (\ref{rad1}) (also see (\ref{rad001})) leads to
$$m-\left(3\lam \mu \int_0^ru(s)^3\;ds\right)^{1/3}r\le u(r)\le m-\left(3\lam \nu(\vep) \int_0^{\vep} u(s)^3\;ds\right)^{1/3}(r-\vep).$$
Hence, $u>\dl$ in some subinterval $[0,t]\subset [0,R]$, where $t>0$. Set
\eqRef{rad4}
r_{\lam}=\sup\{r:\;u(t)>\dl,\;\;0\le t<r\le R\}.
\ee

\NI{\bf Step 4.} From (\ref{rad1}) and Step 3, it is clear that $u\in C(\overline{B_{r_{\lam}}(o)})$ and solves
\eqRef{rad5}
\Df u+\lam a(x) u^3=0,\;\mbox{in $B_{r_{\lam}}(o),$ and $u\ge \dl$ on $\p B_{r_{\dl}}(o)$.}
\ee
We also note that any positive scalar multiple of $u$ also solves (\ref{rad1}). For Cases 1 and 2, we assume that $\dl>0$ and $\lam<\lam_B.$
\vsp
\NI{\bf Case 1:} If $r_{\lam}=R$ then by (\ref{rad4}), $u(R)\ge \dl$. If $u(R)>\dl$, scale $u$ such that $u(R)=\dl$. This provides us with the unique solution to 
(\ref{rad2}), see Remark \ref{unq}.
\vsp
\NI{\bf Case 2:} Suppose that $r_{\lam}<R$. By the continuity of $u$, $u(r_{\lam})=\dl$. We continue $u$ past $r_{\dl}$, using (\ref{rad1}). If $u(r)>0,\;r_{\lam}<r\le R,$ then we scale $u$ such that $u(R)=\dl$. Suppose that there is an $\bar{r}$ with $r_{\lam}<\bar{r}\le R$ such that $u(\bar{r})=0$ (see the estimate in Step 3).  Then $u>0$ in $B_{\bar{r}}(o)$ and satisfies the differential equation in (\ref{rad5}), in $B_{\bar{r}}(o)$, with $u(\bar{r})=0$. If $\bar{r}<R$, by Remarks \ref{mon} and \ref{eigv1}, we have that $\lam\ge\lam_B$, a contradiction. If $\bar{r}=R$, then $u=0$ in $B_{\bar{r}}(o)$, by Lemma \ref{nst}. Thus
$u>0$ in $0\le r\le R$. We may now scale $u$ such that $u(R)=\dl.$ Uniqueness follows from Remark \ref{unq}. This proves part (i).
\vsp
\NI{\bf Step 5.}  Fix $\dl>0$. For each $0<\lam<\lam_B$, part (i) provides us with a unique solution to (\ref{rad2}) which we label as $u_{\lam}$. The function $u_{\lam}$ is positive and radial. 
As has been shown, $u_{\lam}$ also solves (\ref{rad1}). Observe that $\sup_{B}u_{\lam}=u_{\lam}(0).$ Working with the functions $v_{\lam}=u_{\lam}/u_{\lam}(0),$ and arguing as in Theorem \ref{eigv}, there is a subsequence $v_{\lam_k}\rightarrow v$, as $\lam_k\rightarrow \lam_B$,  where $v$ is in $C(\overline{B})$ and solves (\ref{rad2}) with $v|_{\p B}=0$. Moreover, by (\ref{sec43}), $\dl/u_{\lam_k}(o)\rightarrow 0$. It is clear that $v$ solves (\ref{rad1}), in $B$ with $m=1$, that is, 
$$v(x)=v(r)=1-(3\lam_B)^{1/3}\int_0^r \left[ \int_0^t a(s) v(s)^3\;ds\right]^{1/3}\;dt,\;\;\mbox{and $v(R)=0$}.$$
Thus $v$ is a first eigenfunction in $B_R(o)$.  Next, if for some $\lam>0$, there is a function $u$, given by (\ref{rad1}) in $B$, that is positive and vanishes on $|x|=R$, then $\lam\ge \lam_B$. This follows from Remark \ref{eigv1} since $u$ solves (\ref{rad2}) with $\dl=0$. Lemma \ref{sgn} now implies that $\lam=\lam_B$.

Part (iii) of the theorem can be obtained by a differentiation. Also see the proof of Lemma \ref{dst}.
$\Box$
\vsp
From hereon we will take $a(r)=1$. Our goal will be to show that, on a ball, an eigenfunction, corresponding to the first eigenvalue, has one sign and all radial solutions are scalar multiples of each other. Let us set
\eqRef{ball0}
\beta=\left(\int^1_{0}\frac{ds}{(1-s^4)^{1/4}}\right)^4.
\ee
\vsp
\begin{rem}\label{ball} In the statement of Theorem (\ref{rad})(part (iii)), if we take $a(x)=1$, $\dl=0$ and $\lam=\lam_B$, then we obtain
$$ \beta=\lam_B R^4.$$
We argue its validity as follows. Take $\lam<\lam_B$. For $\dl>0$, the corresponding solution $u$ is positive and unique, see Lemma \ref{unq}. We now recall the argument used  in Theorem \ref{eigv}. Set $m=\sup_B u$ and recall that $m=m(\dl)$ becomes unbounded as $\lam\rightarrow \lam_B$, see (\ref{sec43}). Thus, $\dl/m\rightarrow 0$ as $\lam\rightarrow \lam_B$. As a matter of fact, $\dl/m$ depends only on $\lam.$ Taking limits in the formula given in part (iii) of Theorem \ref{rad}, the formula for $\lam_B$ holds.
Also the eigenfunction $u(x)=u(r),\;r=|x|$, given by Theorem \ref{rad}, satisfies the radial version of (\ref{rad2}), that is,
\eqRef{ball1}
\left(\frac{du}{dr}\right)^2\frac{d^2u}{dr^2}+\frac{\beta}{R^4}u^3=0,\;\mbox{in $B_R(o)$, $u^{\prime}(0)=0$ and $u(R)=0$.}\;\;\;\;\Box
\ee
\end{rem}
\vsp
 We now show that the eigenvalue problem on the ball has infinitely many eigenvalues. We also compute the first eigenvalue of an annulus.
For $0\le \kappa<\tau<\infty$ and $p\in \IR^n$, let $\Om=B_{\tau}(p)\setminus \overline{B_{\kappa}(p)}$ be the spherical annulus centered at $p$. Set $2\rh=\tau-\kappa$ and $B=B_{\rh}(p).$ One of our results shows that $\lamo=\lam_B=\beta \rho^{-4}.$ Since, $\Om$ contains a ball of the same size as $B$, this shows that there is no strict domain monotonicity, in general. We refer the reader to Remark \ref{mon}. 

\begin{lem}\label{annb}
Let $R>0$, $p\in \IR^n$ and $\beta$ be as in (\ref{ball0}). Then the problem 
$$\Df u+\lam u^3=0,\;\;\mbox{in $B_R(p)$, and $u=0$ on $\p B_R(p)$,}$$
has infinitely many eigenvalues $\lam$. Moreover, the following hold.\\
\NI (i) The eigenvalues given by
$\lam_{\ell}= \beta (2\ell-1)^4 R^{-4},\;\;\ell=1,2,\cdots,$
have corresponding radial eigenfunctions. \\ 
\NI (ii) Let $0\le \kappa<\tau<\infty$, $p\in \IR^n$, and $\Om=B_{\tau}(p)\setminus \overline{B_{\kappa}(p)}$ be the spherical annulus centered at $p$. Set $2R=\tau-\kappa$ and $B=B_{R}(p).$ Then
$\lamo=\lam_B=\beta R^{-4}.$\quad $\Box$ 
\end{lem}

\NI{\bf Proof.} We carry out the proof in four steps. We refer to Theorem \ref{rad} for the existence of a radial first eigenfunction, also see (\ref{ball1}). The proof has ideas similar to those in Lemma 6.1 in \cite{BMO2}. Set $u^{\prime}(r)=du/dr.$

\NI{\bf Step 1:} Set $r=|x-p|$, and let $u(x)=u(r),\;0\le r\le R$, be a positive radial first eigenfunction of $\Df$ on $B$.  We scale $\sup_B u=1$, and extend $u$ to the rest of $\IR^n$ as follows.  To aid our construction, we recall (\ref{rad1}) and set $a(x)=1$, that is, 
\eqRef{ann01}
u(x)=u(r)=1-(3\lam_B)^{1/3}\int_0^r \left[ \int_0^t u(s)^3\;ds\right]^{1/3}\;dt,\;\;\mbox{and $u^{\prime}(0)=u(R)=0.$}
\ee
First we use an odd reflection about $r=R$. Define
$$u_1(r)=\left\{\begin{array}{ccl} u(r),&\;0\le r\le R,\\ -u(2R-r),&\; R\le r\le 2 R.\end{array}\right.$$
Thus, $u_1$ satisfies
$$\left(\frac{du_1}{dr}\right)^2\frac{d^2u_1}{dr^2}+\lam_B u_1^3=0,$$
in $(0,2R)$, except perhaps at $r=R.$ Next we use an even reflection about $r=2R$ and define
$$u_2(r)=\left\{\begin{array}{ccl} u_1(r),&\;0\le r \le 2R,\\ u_1(4R-r),&\;2R\le r\le 4R.\end{array}\right.$$
Finally, we use a $4R$-periodic extension of $u_2$ to all of $[0,\infty)$. More precisely, for $0\le r<\infty$, let $k=1,2,\cdots,$ be such that $4kR \le  r \le 4(k+1)R$. Now, define
$$u_{\infty}(r)=u_2(r-4kR),\;\;\; \mbox{for}\;4kR\le r\le 4(k+1)R.$$
\NI{\bf Step 2:} Our goal is to show that $u_{\infty}$ solves 
\eqRef{ann0}
\Df u_{\infty}+\lam_B u_{\infty}^3=0,\;\;\mbox{in $\IR^n$.}
\ee
It is clear from Step 1 that we need check this assertion only at $r=R,\;2R$.  We prove this first for $r=R$. We work with $u_1$. Suppose that $\psi\in C^2(\IR^n)$, and $u_1-\psi$ has a maximum at a point $q\in \p B_{R}(p).$ 
We may assume that the segment $pq$ lies along the positive $x_n$ axis. Let $e_n$ denote the unit vector along the positive $x_n$ axis. 
By our construction and (\ref{ann01}), $u_1(q)=u_1(R)=0$, thus implying that 
\eqRef{ann00}
u_1(x)\le \psi(x)-\psi(q)=\langle D\psi(q), \;x-q\rangle+o(|x-q|),\;\;\;\mbox{as}\;x\rightarrow q.
\end{equation}
Take $x\in \p B_{R}(p)$. Since $u_1(x)=0$, dividing both sides by $|x-q|$ and letting $x\rightarrow q$, we get $D\psi(q)=\pm |D\psi(q)|e_n.$ Next,
for small $\tht$, select $x=q+\tht e_n$. Since $u(r)=u(R+\tht)$, we have
$$\frac{u_1(R+\tht)}{|\tht|}\le \frac{\tht}{|\tht|}\langle D\psi(q),e_n\rangle+o(1),\;\;\mbox{as}\;\tht\rightarrow 0.$$
We select $\tht<0$ and note that $u_1(R+\tht)>0$, see (\ref{ann01}) and Step 1. Noting that $u_1^{\prime}(R)=u^{\prime}(R-)$, we get $\langle D\psi(p),e_n\rangle \le u_1^{\prime}(R).$ Now choosing $\tht>0$ and recalling that $u_1(R+\tht)<0$, we obtain $D\psi(q)=u_1^{\prime}(R)e_n$. 
Next, a simple calculation leads to $\Df \psi(p)=(u_1^{\prime}(R))^2D_{nn}\psi(p).$ To determine the sign of $D_{nn}\psi(p)$, we use (\ref{ann00}) to obtain
\eqRef{ann1}
u_1(x)\le \langle u_1^{\prime}(\rh)e_n, x-q\rangle +\frac{\langle D^2\psi(q)(x-q), x-q\rangle}{2}+o(|x-q|^2),\;\;\;\mbox{as}\;x\rightarrow q.
\ee
Taking $x=q+\tht e_n$, where $\tht$ is small, it follows that
$$u_1(R+\tht)\le \tht u_1^{\prime}(R)+\left(\frac{\tht^2}{2}\right)D_{nn}\psi(q)+o(\tht^2),\;\;\mbox{as}\;\tht\rightarrow 0.$$ 
Using (\ref{ann01}) and Step 1, a differentiation yields that $u_1^{\prime\prime}(R-)=u_1^{\prime\prime}(R+)=0.$ Clearly, $u_1$ is $C^2$ near $r=R$, if we define
$u_1^{\prime\prime}(R)=0$. 
Using Taylor's expansion of $u_1$ at $r=R$, we obtain, for small $\tht$, 
$$u_1(R+\tht)-\tht u_1^{\prime}(R)=\frac{\tht^2}{2}u^{\prime\prime}(R)+o(\tht^2)\le \frac{\tht^2}{2} D_{nn}\psi(q)+o(\tht^2),\;\;\;\mbox{as}\;\tht\rightarrow 0.$$ 
Hence, $D_{nn}\psi(q)\ge 0$ and now recalling that $u(q)=0$, we have  $\Df\psi(q)+\lam_B u^3(q)=(u_1^{\prime}(R))^2D_{nn}\psi(q)\ge 0.$ Thus $u_{\infty}$ is a sub-solution near $|x|=R.$ 

Now suppose that for some $\psi\in C^2$, $u_1-\psi$ has a minimum at some $q\in \p B_{R}(p)$. Then $(-u_1)-(-\psi)$ has a maximum at $q$. Arguing as above we conclude that $u_1$ is a super-solution near $|x|=R$. 

 To prove that $u_{\infty}$ solves (\ref{ann0}) near $|x|=2 R$, we observe that $u_{\infty}^{\prime}(0+)=u_{\infty}^{\prime}(2 R)=0.$
This together with the arguments employed in Theorem \ref{rad}(see Step 2) may be now used to treat the case $r=2 R.$ Thus (\ref{ann0}) holds. 
\vsp
\NI{\bf Step 3.} From our construction of $u_{\infty}$ in Step 1, it is clear that $u_{\infty}((2\ell-1)R)=0$, for $\ell=1,2,\cdots.$ Next, by a differentiation, we see that 
the function $w(r)=u_{\infty}((2\ell-1)r)$ provides us with an eigenfunction on $B_{R}(p)$ 
corresponding to the eigenvalue $\lam_{\ell}=(2\ell-1)^4\beta/R^4.$ 
This proves part (i).
\vsp
\NI{\bf Step 4:} We now address part (ii) of the lemma. Recall that $\Om=\{x:\;\kappa<|x|<\tau\}$ and $2R=\tau-\kappa$. If, for some $\ell=0,1,2,\cdots$, $\kappa=(2\ell+1)R$ then 
$\tau=(2\ell+3)R$. From Step 1, for every $\ell$, $u_{\infty}((2\ell+1)R)=0$ and $u_{\infty}$ has one sign in $[(2\ell+1)R,\;(2\ell+3)R]$. Hence, $u_{\infty}=0$, on $\p\Om$, and $u_{\infty}$ has one sign in $\Om$. Thus, (\ref{ann0}) and Lemma \ref{sgn} imply that $u_{\infty}(r)$, restricted to $[\kappa,\;\tau]$, is a first eigenfunction on $\Om$, and $\lamo=\lam_B$. If $(2\ell+1)R<\kappa<(2\ell+3)R$ for some $\ell$, then the function 
$v(r)=u_{\infty}(r-\dl),\;\dl=\kappa-(2\ell+1)R$ is a first eigenfunction in $[\kappa,\;\tau]$. 
If $0\le \kappa<R$ and $\dl=R-\kappa$, then $v(r)=u_{\infty}(r+\dl),\;\kappa\le r\le \tau,$ is the desired eigenfunction.
Note that the ordinary differential equation in Remark \ref{ball1} is translation invariant. In any case, $\lamo=\lam_B$.  Note also that if $A=B_{2\rh}(p)\setminus\{p\}$ then $\lam_A=\lam_B$.  $\qed$ 
\vsp
Finally, we prove that the first eigenfunction, on the ball, has one sign and that all a radial solution is unique up to scalar multiplication. Simplicity of $\lam_B$ would follow if every solution is radial. However,  it is not clear to us if this is indeed true.
\begin{thm}\label{simp}
Let $R>0$, let $u\in C(\overline{B_R(o)})$ solve the eigenvalue problem
$$\Df u+\lam_B u^3=0,\;\mbox{in $B_R(o)$, and $u=0$ on $\p B_R(o)$.}$$
It follows that (i) $u$ has one sign in $B_R(o)$, and (ii) if $u$ is radial and $\sup_B u=1$ then $u$ is unique. 
\end{thm}

\NI{\bf Proof.} Set $B=B_R(o)$; scale $u$ so that $u(o)=1$.  Set $B^+=\{x\in B:\;u(x)>0\}$ and $B^-=\{x\in B:\;u(x)<0\}$. Note that $u$ is infinity super-harmonic in $B^+$ and infinity sub-harmonic in $B^-$.

We prove part (i). Assume that $u$ changes sign in $B$. We discuss the case when a component $C$ of $B^-$ is compactly contained in $B$. Since $u$ is an eigenfunction on $B^-$, Remark \ref{eigv1} implies that
$\lam_{C}\le\lam_B$. This contradicts the strict monotonicity shown in Remark \ref{mon}.  Thus, if $B^-$ is non-empty then $B^-\cap \p B_r(o)$ is non-empty, for every $r$ close to $R$. 

We derive bounds for $u$. Set $m=\inf_B u$ and $M=\sup_B u$. By our hypothesis, $m<0<M$. For any $L\not=0$, select $b=b(L)$ such that 
$(b+3LR)^{4/3}-b^{4/3}=-4L.$
Then the function 
$$\psi(x)=\psi(|x|,L)=1+\frac{1}{4L}\left[\left(b+3L|x|\right)^{4/3}-b^{4/3}\right],\;\;\mbox{for}\;x\in B,$$
satisfies 
$$\Df \psi=L,\;\mbox{in $B\setminus\{o\}$, $\psi(0)=1$ and $\psi(R)=0.$}$$   
Set $\psi_M(x)=\psi(|x|,-8\lam_BM^3)$ and $\psi_m(x)=\psi(|x|, 8\lam_B |m|^3)$. Then 
$$\Df \psi_M=-8\lam_B M^3\le \Df u=-\lam_B u^3\le \Df \psi_m=8\lam_B |m|^3.$$ 
Since $u(o)=\psi_m(0)=\psi_M(0)=1$, and $\psi_m(R)=\psi_M(R)=u=0$ on $\p B$, Lemma \ref{cmp3} implies
\eqRef{simp1}
\psi_m(x)\le u(x)\le \psi_M(x),\;\;\mbox{for $x\in \overline{B}$}.
\ee
Consider all rotations of $B$ about $o$. Let $A$ be an $n\times n$, orthogonal matrix. Define
$$u_{\ell}(x)=\inf_{A} u(Ax),\;\;\mbox{for $x\in \Om$.}$$ 
Set $r=|x|$, clearly, $u_{\ell}(x)=u_{\ell}(r)=\inf_{\p B_r(o)} u$ and $u(o)=u_{\ell}(o)=1$. Since $\Df$ is rotation invariant, $u(Ax)$ is an eigenfunction. Arguing as in Theorem 3.1 in \cite{BMO2} (this appears in the Perron method and uses a perturbation, see equations (3.2)-(3.4) therein),  $u_{\ell}$ is a super-solution, that is, 
\ben
\Df u_{\ell}+\lam_B u_{\ell}^3\le 0,\;\mbox{in $B$, and $u_{\ell}(R)=0.$}
\een
Since every $u(Ax)$ satisfies (\ref{simp1}), we have $\psi_m(x)\le u_{\ell}(x)\le \psi_M(x)$ in $\Om$. Thus, by Lemma \ref{lpc}, $u_{\ell}$ is locally Lipschitz continuous in $B$, and
$u_{\ell}(r)$ assumes the zero boundary data continuously.

Define $r_0=\sup\{r:\;u_{\ell}(t)>0,\;\forall\;0\le t<r\}$. Recalling that $u(o)=1$ and $B^-$ is non-empty, we see that $0<r_0<R$ and $u_{\ell}(r_0)=0.$ Since every component of $B^-$ meets $\p B$, $u_{\ell}(r)<0$ in $r_0<r<R$. Set $A=\{x:\;r_0<|x|<R\}$ and $d=(R-r_0)/2$. 
We take $v=-u_{\ell}$ to obtain
$$\Df v+\lam_B v^3\ge 0,\;v>0,\;\mbox{in $A$ and $v=0$ on $\p A.$}$$
Next, Lemma \ref{annb} implies that $\lam_A=\lam_{B_d(o)}>\lam_B$.  By Remark \ref{set21} and Lemma \ref{nst}, 
we get $v\le 0$ in $A$. This is a contradiction and it follows that $u\ge 0$ in $B$, and hence, $u>0$. 

We now prove part (ii). Let $v>0$ be the radial solution in $B$ given by part (ii) of Theorem \ref{rad}.  Suppose that $u$ is a radial first eigenfunction on $B$.
Since every eigenfunction has one sign, we may use Remark \ref{unq}  in $B_r(o),\;0<r<R$. Thus
$$\inf_{\p B_r(o)}\frac{u}{v}=\frac{u(r)}{v(r)}\le \frac{u(t)}{v(t)}\le \frac{u(r)}{v(r)}=\sup_{\p B_r(o)}\frac{u}{v},\;\;0\le t\le r.$$ 
Thus $v(r)=(v(o)/u(o))u(r)$ for any $0\le r<R$. Thus $u$ is a scalar multiple of $v$ and uniqueness follows. $\Box$

\vsp
\NI Department of Mathematics\\
\NI Western Kentucky University\\
\NI Bowling Green, Ky 42101
\vsp
\NI Department of Mathematics\\
\NI University of Kentucky\\
\NI Lexington, KY 40506-0027

\end{document}